\documentclass[a4paper,12pt]{article}
\usepackage{amsfonts}
\usepackage{amsmath}
\usepackage{CJK,graphicx}
\usepackage{amscd}
\usepackage{amssymb}
\begin{document}
\title{\textbf{\Huge{The Gauss-Bonnet-Chern Theorem on Riemannian Manifolds}}}\author{Yin Li}\date{}\maketitle
\begin{abstract}
{This expository paper contains a detailed introduction to some important works concerning the Gauss-Bonnet-Chern theorem. The study of this theorem has a long history dating back to Gauss's Theorema Egregium (Latin: Remarkable Theorem) and culminated in Chern's groundbreaking work \cite {sc1} in 1944, which is a deep and wonderful application of Elie Cartan's formalism. The idea and tools in \cite{sc1} have a great generalization and continue to produce important results till today.\\
In this paper, we give four different proofs of the Gauss-Bonnet-Chern theorem on Riemannian manifolds, namely Chern's simple intrinsic proof, a topological proof, Mathai-Quillen's Thom form proof and McKean-Singer-Patodi's heat equation proof. These proofs are related with remarkable developments in differential geometry such as the Chern-Weil theory, theory of characteristic classes, Mathai-Quillen's formalism and the Atiyah-Singer index theorem. It is through these brilliant achievements the great importance and influence of Chern's insights and ideas are shown. Our purpose here is to use the Gauss-Bonnet-Chern theorem as a guide to expose the reader to some advanced topics in modern differential geometry.}
\end{abstract}
\newpage
\begin{center}
\textbf{\Huge{In Memory of Professor S.S. Chern (1911-2004)}}\\

\Large{\textmd{The master does his job and then stops.\\
He understands that the universe is forever out of control,\\
and that trying to dominate events goes against the current of the Tao.\\
Because he believes in himself,\\
he doesn't try to convince others.\\
Because he is content with himself,\\
he doesn't need others' approval.\\
Because he accepts himself,\\
the whole world accepts him.}}
\end{center}
\newpage
\tableofcontents
\section{Introduction}

\subsection{Two Viewpoints of the Gauss-Bonnet Theorem}
Let $M$ be a closed oriented Riemannian surface and $K$ its Gaussian curvature, $P: F\rightarrow M$ a diffeomorphism of a polygon $F$ onto a subset of $M$, $\alpha_i$ the exterior angles of the vertices of $P(F)$, and $\kappa_g$ the geodesic curvature of the positively oriented curve $\partial P$. The classical Gauss-Bonnet theorem says that (\cite{ccl}):
\begin{equation}\label{eq:ori}\int_PK\textrm{d}A+\int_{\partial P}\kappa_g\textrm{d}s+\sum_i\alpha_i=2\pi.\end{equation}
Given a triangulation of $M$, it is then easy to deduce from (\ref{eq:ori}) that
\begin{equation}\label{eq:ori1}\int_MK\textrm{d}A=\chi(M),\end{equation}
where $\chi(M)$ denotes the Euler characteristic of $M$.\\
(\ref{eq:ori}) and (\ref{eq:ori1}) respectively give the two viewpoints of the Gauss-Bonnet theorem in history.\\
In fact, let $P$ be a geodesic triangle, then (\ref{eq:ori}) becomes
\begin{equation}\label{eq:aw1}\int_PK\textrm{d}A+\sum_{i=1}^3\alpha_i=2\pi,\end{equation}
so it reveals the relationship between the Gaussian curvature and the ``angular excess", i.e., the sum of the interior angles of a geodesic triangle minus the sum of the interior angles of a triangle in $\mathbb{R}^2$.\\
On the other hand, (2) gives the relatively modern viewpoint of the Gauss-Bonnet theorem, which says that the most important topological invariant of $M$, i.e., the Euler characteristic, is given by the so called ``curvature integral". For the fundamental role the Euler characteristic plays in topology and geometry, I strongly recommend the reader to read \cite{sc3}.\\
Heinz Hopf proposed in the late 1920s the question of generalizing the Gauss-Bonnet theorem (\ref{eq:ori1}) to all even dimensions. This problem was first solved by Allendoerfer and Weil in 1943. In \cite{aw}, they proved the following analogue of (\ref{eq:ori}) for all Riemannian manifolds:
\begin{equation}\label{eq:awf}(-1)^d\chi'(P)=\int_P\Psi(z)\textrm{dvol}_z+\sum_i\int_{\partial P_i}\int_{\mathit{\Gamma}(\zeta)}\Psi(\zeta,d\xi|\partial P_i)\textrm{dvol}_\zeta,\end{equation}
where
\begin{eqnarray}
\Psi(z)=(2\pi)^{-\frac{d}{2}}\frac{1}{2^d(\frac{d}{2})!g}\cdot\qquad\qquad\qquad\qquad\qquad\qquad\qquad\qquad\qquad\quad\nonumber\\
\sum_{\sigma_1,\sigma_2\in\Sigma_d}\textrm{sgn }\sigma_1\cdot\textrm{sgn }\sigma_2R_{\sigma_1(1)\sigma_1(2)\sigma_2(1)\sigma_2(2)}\cdot\cdot\cdot R_{\sigma_1(d-1)\sigma_1(d)\sigma_2(d-1)\sigma_2(d)},
\end{eqnarray}
$d=\dim_{\mathbb{R}}M$, $R_{ijkl}$ is the Riemann curvature tensor with respect to the metric $g_{ij}$, $g=\det(g_{ij})$. $\Sigma_d$ denotes the set of all the permutations of $\{1,\cdot\cdot\cdot,d\}$ and sgn denotes the signature of a given permutation. dvol is the Riemannian volume element, $P$ is a differentiable Riemannian polyhedron, whose boundary consists of Riemannian polyhedra $\partial P_i$, $\mathit{\Gamma}(\zeta)$ is a subset of the unit sphere centering at $\zeta$ which generalizes the exterior angle of $\partial P_i$ at $\zeta$ in the case of a Riemannian surface, the term $\Psi(\zeta,d\xi|\partial P_i)$ is a measure of the curvature of $\partial P_i$ which generalizes the geodesic curvature of a curve on a Riemannian surface, and $\chi'$ is the interior Euler characteristic, which is the Euler characteristic computed only on the interior faces of a triangulation.\\
Actually, (\ref{eq:aw1}) is the Gauss-Bonnet theorem for Riemannian manifolds with boundary. Using (\ref{eq:aw1}), they then deduced the Gauss-Bonnet theorem for Riemannian manifolds without boundary as follows, which is a generalization of (\ref{eq:ori1}):
\begin{equation}\label{eq:gbc}\int_M\Psi(z)\textrm{dvol}_z=\chi(M).\end{equation}
However, the paper of Allendoerfer and Weil is very complicated and involves equally complicated works of Allendoerfer (\cite{ca}), Fenchel (\cite{wf}) and Weyl (\cite{hw}). Moreover, they used in their proof the local isometric embedding of the Riemannian manifolds into the Eucildean spaces, which seems unnatural for proving an intrinsic formula.\\
It was Chern who made the Gauss-Bonnet theorem (\ref{eq:gbc}) widely known. His simple intrinsic proof in 1944 not only gave a beautiful and profound proof of the Gauss-Bonnet theorem, but also enlightened the whole field of global differential geometry. In his 6-page paper \cite{sc1}, he invented the technique now called transgression, introduced curvature to topology, and showed the importance and usefulness of the concept of a fiber bundle in differential geometry.\\
It is the aim of this expository paper to review Chern's great paper \cite{sc1} together with some other important developments on the Gauss-Bonnet-Chern theorem. As an undergraduate student interested in differential geometry, the author write this paper to learn, and also to dedicate it to the memory of Professor S.S. Chern on the occasion of his 100th birthday. I sincerely hope that he would like it.
\subsection{Arrangement of This Paper}
I have endeavored to write this article in a self-contained way, only assume the reader know elementary analysis, geometry and topology. The reader who don't familiar with these materials may refer to \cite{bt}, \cite{jj} and \cite{jm}. We shall work with real coefficients unless otherwise mentioned.\\
For clarity, we state two basic topological facts in $\S2$, namely the elementary singularity theory of a unit vector field and the Poincar\'{e}-Hopf index theorem. They will be used in Chern's intrinsic proof in $\S3$.\\
$\S3$ is devoted to Chern's paper \cite{sc1}. However, we will briefly mention the Chern-Weil theory at the end of this section since it is a natural generalization of the Gauss-Bonnet-Chern theorem under the general framework of the Cartan formalism and gives enough motivations for the topological approach to the Gauss-Bonnet-Chern theorem which we will discuss in $\S4$. A brief account of Cartan's method of moving frames is given at the beginning of this section since most of the text books on Riemannian geometry are written in the language of tensor analysis.\\
$\S4$ contains two proofs of the Gauss-Bonnet-Chern theorem using the Euler form and the Thom form respectively, in the spirit of the Chern-Weil theory. The first proof is a new proof based on an explicit expression of the Euler class by transition functions and generalizes the Gauss-Bonnet-Chern theorem for a metric-compatible connection on oriented Riemannian vector bundles. The second proof follows directly from Mathai and Quillen's geometric construction of the Thom form in \cite{mq}. We shall develop in detail the topological backgrounds for understanding these proofs at the beginning of $\S4$.\\
In $\S5$ we present an analytic approach to the Gauss-Bonnet-Chern theorem, namely by using the heat equations. Although this proof is a little lengthy, it is a powerful method in mordern differential geometry. We will construct the parametrix and the heat kernel in full detail and the McKean-Singer conjecture will be raised in a natural way. The proof will be completed by using a classical tensor calculus method developed by Patodi in his phD thesis \cite{vp}, since this gives the most self-contained proof.\\
Although this paper contains no new result, several proofs are new and some simplifications have been made to original papers. Noteworthy is our treatment for the Euler class in $\S4.1$, which is different from other standard references, it simplifies our exposition in many places.
\section{Topological Preliminaries}
References for this section are \cite{agp}, \cite{bt}, \cite{mt} and \cite{jm}.
\subsection{Singularities of A Unit Vector Field}
Let $S$ be an oriented sphere bundle over a differentiable manifold $M$. Although in general $S$ does not admit a global section, a well-known result in differential topology says that there may be a section $s$ over the complement of an isolated set $I\subset M$, i.e., $s\in\Gamma(M\setminus I,S)$. The points in $I$ are called singularities of $s$.\\
\textbf{Theorem 2.1.1.} \textit{Let} $\pi:S\rightarrow M$ \textit{be a} $(d-1)$-\textit{sphere bundle over a closed manifold} $M$ \textit{of dimension} $d$. S\textit{uppose that the structure group of} $S$ \textit{can be reduced to the orthogonal group} $\textrm{O}(d,\mathbb{R})$, \textit{then there exists a smooth map} $s:M\rightarrow S$ \textit{such that} $s\in\Gamma(M\setminus I,S)$, \textit{where} $I$ \textit{is discrete}.\\
\textsc{Remarks}
\begin{itemize}
\item Under the assumption of Theorem 2.1.1, the sphere bundle $S$ is induced by a vector bundle $E$. Suppose $s'\in\Gamma(M\setminus I,S)$, then we can set $s=s'$ on $M\setminus I$ and $s=0$ on $I$ to produce an $s\in\Gamma(M,E)$ with the zero set $I$. Thus Theorem 2.1.1 can be restated for vector bundles and zeros of its sections.
\item It follows from the basic facts of the obstruction theory that the assumption ``the structure group of $S$ can be reduced to $\textrm{O}(d,\mathbb{R})$" can be removed. See \cite{ns}.
\end{itemize}
\subsection{Poincar\'{e}-Hopf Index Theorem}
The Poincar\'{e}-Hopf index theorem plays a pivotal role in Chern's proof of the Gauss-Bonnet theorem in \cite{sc1}, since it localizes the global topological information of the manifold using the zeros of a vector field. To state this theorem, we need some basic concepts which give analytic descriptions of the zeros of a vector field.\\
\textbf{Definition 2.2.1.} Let $f:M\rightarrow N$ be a smooth map between two closed oriented manifolds of dimension $d$. Then $f$ defines a pullback on de Rham cohomology $f^\ast:H_{\textrm{dR}}^d(N,\mathbb{R})\rightarrow H_{\textrm{dR}}^d(M,\mathbb{R})$. Let $\omega$ be the generator of $H_{\textrm{dR}}^d(N,\mathbb{R})$, then the \textit{mapping degree} of $f$ is defined to be $\int_Mf^\ast\omega$.\\
Suppose $S$ is a sphere bundle over a closed oriented $d$-dimensional Riemannian manifold $M$ and $s\in\Gamma\Big(U\setminus\{x\},S\Big)$. Suppose $U$ is chosen small enough so that $U\setminus\{x\}\simeq\mathbb{B}^d\setminus\{0\}$ and $S|_{U\setminus\{x\}}\simeq \Big(U\setminus\{x\}\Big)\times\mathbb{S}^{d-1}$, where $\simeq$ denotes the diffeomorphism and $\mathbb{B}^d\subset\mathbb{R}^d$ is the standard Euclidean ball. Let $B_r\subset M$ denote the preimage of $r\mathbb{B}^d$ under the trivialization of $M$, then $B_r$ is clearly oriented. Choose the orientation on $\mathbb{S}^{d-1}$ such that the diffeomorhism $S|_{B_r}\simeq B_r\times\mathbb{S}^{d-1}$ is orientation preserving, where $B_r\times\mathbb{S}^{d-1}$ is given the product orientation.\\
\textbf{Definition 2.2.2.} The \textit{local degree} of the section $s\in\Gamma\Big(U\setminus\{x\},S\Big)$ at $x$ is defined to be the mapping degree of the composition
\begin{equation}\partial\overline{B_r}\xrightarrow{s}S|_{B_r}=\overline{B_r}\times\mathbb{S}^{d-1}\xrightarrow{p}\mathbb{S}^{d-1},\end{equation}
where $p$ is the projection and $\overline{B_r}$ is the closure of $B_r$.\\
\textbf{Definition 2.2.3.} Let $X$ be a vector field with discrete zeros on $M$, then the \textit{index} of $X$ at a zero $x\in I$ is defined to be the local degree at $x$ of the section $\frac{X}{\|X\|}\in\Gamma(M\setminus I,SM)$, where $SM$ denotes the projective sphere bundle relative to some Riemannian metric of $M$.\\
Denote by $\textrm{ind}_X(x)$ the index of the vector field $X$ at $x$, we state the Poincar\'{e}-Hopf index theorem as follows:\\
\textbf{Theorem 2.2.4. (Poincar\'{e}-Hopf)} \textit{Let the set} $I$ \textit{be defined as above and suppose it is chosen to be discrete}. $M$ \textit{is a closed oriented manifold and} $X$ \textit{is a vector field on} $M$, \textit{then}
\begin{equation}\label{eq:phi}\chi(M)=\sum_{x_i\in I}\textrm{ind}_X(x_i).\end{equation}
\textsc{Remarks}
\begin{itemize}
\item An analytic proof of this theorem was obtained by Witten in his famous paper \cite{ew}, see also \cite{wz} for an exposition of his proof. An alternative analytic proof based on the idea of \cite{ma} was found in \cite{wz2}.
\item The Poincar\'{e}-Hopf index theorem is the simplest example of the localization theorems which relate the characteristic numbers with zeros of a vector field. The most famous example among these theorems is the Bott residue formula \cite{rb}. However, it turns out that the Bott residue formula holds even when the zero set of a vector field is not discrete, this amazing extension was made in \cite{bc}. See also \cite{bgv}, \cite{sc4} and \cite{wz} for expositions of these results.
\item An extension of this theorem will be made in $\S4.2$ as a byproduct of the generalized Gauss-Bonnet-Chern theorem, see (\ref{eq:gph}).
\end{itemize}
\section{A Simple Intrinsic Proof}

\subsection{Maurer-Cartan Equations}
Here we recall briefly the method of moving frames and the corresponding structure equations. Of course the standard references for these topics are \cite{ec} and \cite{ccl}. For a modern introduction to differential geometry via moving frames, the reader may refer to \cite{il} and \cite{rs}.\\
We begin with the general case when $\pi:P\rightarrow M$ is a principal $G$-bundle over a $d$-dimensional differentiable manifold $M$. Let $\mathfrak{g}$ denote the Lie algebra of $G$, then it can be identified with the tangent space $T_eG$, where $e$ is the unit element of $G$. Suppose $\{U_\alpha\}$ is an open cover of $M$ on which $P$ is trivial, then there exist diffeomorphsims $\varphi_\alpha:\pi^{-1}(U_\alpha)\rightarrow U_\alpha\times G$ such that
\begin{equation}\label{eq:triv}\varphi_\alpha(y)=\Big(\pi(y),p_\alpha(y)\Big),p_\alpha(g\cdot y)=g\cdot p_\alpha(y),\forall g\in G.\end{equation}
The group formed by the transition functions $\{g_{\alpha\beta}\}$ of $P$ can be identified with a subgroup of $G$, then we have in $G$ the relation
\begin{equation}\label{eq:pg}p_\alpha(y)g_{\alpha\beta}=p_\beta(y).\end{equation}
Let $R_g$ be the right translation $h\rightarrow h\cdot g$, where $g,h\in G$, and $\textrm{Ad}(g)$ be the adjoint action in $\mathfrak{g}$ defined by
\begin{equation}\exp_e\Big(\textrm{Ad}(g)X\Big)=g(\exp_eX)g^{-1},g\in G,X\in\mathfrak{g},\end{equation}
we now introduce the concept of a Cartan connection.\\
\textbf{Definition 3.1.1.} Let $\pi:P\rightarrow M$ be a principal $G$-bundle, a \textit{Cartan connection} $\nabla^P$ on $P$ is given by a $\mathfrak{g}$-valued 1-form $\theta_\alpha$ on every $U_\alpha$ such that
\begin{equation}\label{eq:conn}\theta_\alpha=(R_{g_{\alpha\beta}^{-1}})_\ast\textrm{Id}(T_{g_{\alpha\beta}}G)+\textrm{Ad}(g_{\alpha\beta})\theta_\beta\textrm{ in }U_\alpha\bigcap U_\beta,\end{equation}
where $\textrm{Id}(T_{g_{\alpha\beta}}G)$ denotes the identity endomorphism of $T_{g_{\alpha\beta}}G$. Note that $(R_{g_{\alpha\beta}^{-1}})_\ast\textrm{Id}(T_{g_{\alpha\beta}}G)$ can be identified with a $\mathfrak{g}$-valued 1-form on $U_\alpha\bigcap U_\beta$.\\
Define $\Theta_\alpha\in\Gamma(U_\alpha,\mathit{\Lambda}^2T^\ast M)\otimes\mathfrak{g}$ by
\begin{equation}\label{eq:mc1}\Theta_\alpha=d\theta_\alpha-\frac{1}{2}[\theta_\alpha,\theta_\alpha],\end{equation}
where $[\cdot,\cdot]$ denotes the Lie bracket for $\mathfrak{g}$-valued differential forms defined by
\begin{equation}\label{eq:conv}[X\otimes\omega,Y\otimes\theta]=[X,Y]\otimes(\omega\wedge\theta),X,Y\in\mathfrak{g},\omega,\theta\in\mathit{\Lambda}^\ast T^\ast M.\end{equation}
(\ref{eq:mc1}) is called the \textit{Maurer-Cartan equations} and $\{\theta_\alpha\}$ are called \textit{Maurer-Cartan forms}. $\{\Theta_\alpha\}$ are by definition the \textit{curvature forms of the Cartan geometry}.\\
Now suppose the restriction $P|_{U_\alpha}$ is a Lie group, then one easily obtains (c.f. \cite{cc})
\begin{equation}d\theta_\alpha=\frac{1}{2}[\theta_\alpha,\theta_\alpha].\end{equation}
This is the Maurer-Cartan equation for Lie groups. It shows that the curvature forms of the Cartan geometry are complete local obstructions to $P$ being a Lie group. By (\ref{eq:conn}), we deduce the local transition formulas
\begin{equation}\label{eq:ltf}\Theta_\alpha=\textrm{Ad}(g_{\alpha\beta})\Theta_\beta\textrm{ in }U_\alpha\bigcap U_\beta.\end{equation}
Exterior differentiation of (\ref{eq:mc1}) gives the \textit{Bianchi identities}
\begin{equation}d\Theta_\alpha=[\theta_\alpha,\Theta_\alpha].\end{equation}
We now apply the above method to Riemannian geometry.\\
Let $M$ be a Riemannian manifold with dimension $d$. In this case, the principal $G$-bundle $P$ over $M$ can be produced by moving an orthogonal frame over $M$, i.e., $G=\textrm{O}(d,\mathbb{R})$. It follows that
\begin{equation}\label{eq:liea}\mathfrak{o}(d,\mathbb{R})=\Big\{X\in\mathfrak{gl}(d,\mathbb{R}):\langle Xv,w\rangle+\langle v,Xw\rangle=0\textrm{ for all }v,w\in\mathbb{R}^d\Big\},\end{equation}
where
\begin{equation}\mathfrak{gl}(d,\mathbb{R})=\{X:\mathbb{R}^d\rightarrow\mathbb{R}^d\textrm{ linear}\}\end{equation}
denotes the Lie algebra of $\textrm{GL}(d,\mathbb{R})$.\\
Therefore $\Big\{g_{\alpha\beta}(x)\Big\}$ can be identified with orthogonal matrices and (\ref{eq:conn}) reduces to the matrix equation
\begin{equation}\label{eq:red}\theta_\alpha=dg_{\alpha\beta}g_{\alpha\beta}^{-1}+\textrm{Ad}(g_{\alpha\beta})\theta_\beta.\end{equation}
This shows that in this case, a connection $\nabla^P$ on $P$ is given by a $d\times d$ matrix of 1-forms $\omega=(\omega_i^{\phantom{i}j})$ on each $U_\alpha$, where $\{U_\alpha\}$ is an open covering of $M$ which trivializes $P$. We shall call $\omega$ the \textit{connection matrix} of $\nabla^P$. By (\ref{eq:liea}), we get in addition that $\omega$ is skew-symmetric, i.e.:
\begin{equation}\label{eq:metr}\omega_i^{\phantom{i}j}+\omega_j^{\phantom{j}i}=0.\end{equation}
Using (\ref{eq:conv}), we write the Maurer-Cartan equations in this case as
\begin{equation}\label{eq:mce2}\Omega_i^{\phantom{i}j}=d\omega_i^{\phantom{i}j}-\omega_i^{\phantom{i}k}\wedge\omega_k^{\phantom{k}j},\end{equation}
where $\Omega_i^{\phantom{i}j}$ are entries of the corresponding \textit{curvature matrix} $\Theta=(\Omega_i^{\phantom{i}j})$. Here and in the sequel, we shall adopt the Einstein summation convention.\\
(\ref{eq:metr}) and (\ref{eq:mce2}) are known as \textit{the fundamental lemma of Riemannian geometry}.\\
By (\ref{eq:mc1}) and (\ref{eq:metr}), it follows that the curvature matrix is also skew-symmetric:
\begin{equation}\label{eq:anti}\Omega_i^{\phantom{i}j}+\Omega_j^{\phantom{j}i}=0.\end{equation}
Similarly, the Bianchi identities now becomes
\begin{equation}\label{eq:bian}d\Omega_i^{\phantom{i}j}+\omega_j^{\phantom{j}k}\wedge\Omega_k^{\phantom{k}i}-\omega_i^{\phantom{i}k}\wedge\Omega_k^{\phantom{k}j}=0.\end{equation}
This is the Riemannian geometry via moving frames.\\
\textsc{Remarks}
\begin{itemize}
\item When $G=\textrm{O}(d,\mathbb{R})$, $\Big\{p_\alpha(y)\Big\}$ in (\ref{eq:triv}) can be identified with orthogonal matrices. Denote by $p_\alpha^1(y)$ the matrix formed by the first row of $p_\alpha(y)$, then (\ref{eq:pg}) restricts to the relation
\begin{equation}\label{eq:trs}p_\alpha^1(y)g_{\alpha\beta}=p_\beta^1(y).\end{equation}
It follows that we can construct a vector bundle $E$ with $\Big\{g_{\alpha\beta}(x)\Big\}$ as its transition functions. $E$ is called the \textit{associated vector bundle} of $P$ and $P$ will be referred to as the \textit{orthogonal frame bundle} of $E$.\\
Combine (\ref{eq:red}) and (\ref{eq:trs}) to get
\begin{equation}\label{eq:idc}(dp_\alpha^1+p_\alpha^1\theta_\alpha)g_{\alpha\beta}=dp_\beta^1+p_\beta^1\theta_\beta.\end{equation}
From (\ref{eq:idc}) it follows that every connection $\nabla^P$ on the principal $\textrm{O}(d,\mathbb{R})$-bundle $P$ induces a connection $\nabla^E$ on its associated vector bundle $E$, and $\nabla^E$ is uniquely determined by
\begin{equation}\nabla^Ep_\alpha^1=dp_\alpha^1+p_\alpha^1\theta_\alpha.\end{equation}
Moreover, by (\ref{eq:metr}), $\nabla^E$ must be metric compatible. This restriction on $\nabla^E$ is essential for our proof of the Gauss-Bonnet-Chern theorem, as we will see later.\\
It is a trivial fact that every metric connection $\nabla^E$ on a vector bundle $E$ also induces uniquely a connection on its orthogonal frame bundle $P$, therefore the equivalence between studying the metric connections on a vector bundle and the method of moving frames has been showed.\\
More generally, for $G=\textrm{GL}(d,\mathbb{C})$, every connection $\nabla^P$ on the principal $\textrm{GL}(d,\mathbb{C})$-bundle $P$ can be identified with a unique connection $\nabla^E$ on the associated Hermitian vector bundle $E$, here $\nabla^E$ may not be metric compatible.
\item We now move the orthogonal frames in the tangent bundle $TM$ and still use $\omega=(\omega_i^{\phantom{i}j})$ and $\Theta=(\Omega_i^{\phantom{i}j})$ to denote the corresponding connection and curvature matrices. As we have seen, this is equivalent to studying a metric connection on $TM$. Now let $\{O,e_1,\cdot\cdot\cdot,e_d\}$ be our local frames and $\mathbf{O}$ be the position vector of $O$. The connection of $e_i\in\Gamma(U_\alpha,TM)$ is given by
\begin{equation}\label{eq:ctm}de_i=\omega_i^{\phantom{i}j}e_j.\end{equation}
Since $\mathbf{O}$ can be identified with a vector in the linear space $\textrm{span}\{e_1,\cdot\cdot\cdot,e_d\}$, we can take its exterior differential
\begin{equation}d\mathbf{O}=\omega^ie_i,\omega^i\in\Gamma(U_\alpha,T^\ast M).\end{equation}
It is easy to verify that if
\begin{equation}d\omega^i=\omega^j\wedge\omega_j^{\phantom{j}i},\end{equation}
then the connection $\nabla^{TM}$ induced by $\nabla^P$, where $P$ is the orthogonal frame bundle of $TM$, is torsion free and is therefore the Levi-Civita connection. We will see that the torsion freeness of $\nabla^{TM}$ is not required in the proof of the Gauss-Bonnet-Chern theorem.\\
Now suppose that $\nabla^{TM}$ is the Levi-Civita connection, let $\omega$ and $\Theta$ be the corresponding connection and curvature matrices respectively. We have the following simple relation between the curvature 2-forms and the components of the Riemann curvature tensor:
\begin{equation}\Omega_i^{\phantom{i}j}={R_i^{\phantom{i}j}}_{kl}d\omega^k\wedge d\omega^l.\end{equation}
In fact, the deep relationship between differential forms and algebraic topology makes the curvature forms more convenient than the curvature tensor in global differential geometry, this is formalized in the Chern-Weil theory in the late 1940s. Chern's simple intrinsic proof of the Gauss-Bonnet theorem is actually a perfect example of the Chern-Weil theory which shows that the Euler class can be expressed by the Pfaffian of the curvature forms. These remarks will be made rigorously in $\S3.3$.
\end{itemize}
\subsection{Chern's Original Proof}
Chern's proof is very beautiful and conceptual. We shall first sketch out the basic steps of his proof and then carry out the details step by step. We assume here that $M$ is a $d$-dimensional closed oriented Riemannian manifold, here $d$ is an even number.\\
\textsc{Sketch of the Proof}
\begin{itemize}
\item Let $I$ be as in $\S2.1$, take an $\mathcal{S}\in\Gamma(M\setminus I,SM)$, i.e., a unit vector field with possibly isolated singularities. The existence of such a section is guaranteed by Theorem 2.1.1. We shall identify $\mathcal{S}$ with $\mathcal{S}(M\setminus I)$.
\item Extract small balls $\bigcup_{x_i\in I}\overline{B(x_i)}$ about these singularities, so $\mathcal{S}$ is well-defined on $M\setminus\bigcup_{x_i\in I}\overline{B(x_i)}$.
\item Define the Gauss-Bonnet integrand $\Omega$, i.e., an intrinsic $d$-form formed by the curvature 2-forms. This is a generalization of the form $KdA$ in the surface case, see (\ref{eq:ori}).
\item Transgression. Use the projection $\pi:SM\rightarrow M$ to pull $\Omega$ back to $SM$ and establish the transgression formula $\pi^\ast\Omega=d\Pi$, where $\Pi\in\Gamma(SM,\mathit{\Lambda}^{d-1}T^\ast SM)$.
\item Use $\pi|_{\mathcal{S}}$, which is the restriction of $\pi:SM\rightarrow M$ on $\mathcal{S}$, to pull the integral $\int_{M\setminus\bigcup_{x_i\in I}\overline{B(x_i)}}\Omega$ back to $\mathcal{S}$. Make the radii of $B(x_i)$ tend to 0, apply the Stokes theorem and the Poincar\'{e}-Hopf index theorem to finish the proof.
\end{itemize}
\textsc{The Gauss-Bonnet Integrand}\\
Let $\sigma\in\Sigma_d$, we construct a differential form $\Omega_\alpha\in\Gamma(U_\alpha,\mathit{\Lambda}^dT^\ast M)$ for every $\alpha$ as follows:
\begin{equation}\label{eq:gbi}\Omega_\alpha=\frac{(-1)^{\frac{d}{2}+1}}{2^d\pi^{\frac{d}{2}}(\frac{d}{2})!}\sum_{\sigma\in\Sigma_d}\textrm{sgn }\sigma\cdot\Omega_{\sigma(1)}^{\phantom{\sigma(1)}\sigma(2)}\wedge\cdot\cdot\cdot\wedge\Omega_{\sigma(d-1)}^{\phantom{\sigma(d-1)}\sigma(d)}.\end{equation}
\textbf{Lemma 3.2.1.} \textit{The form} $\Omega_\alpha$ \textit{constructed above is intrinsic and global}.\\
\textit{Proof}. Let $A=(a_{ij})$ be an orthogonal transformation which changes $\{O,e_1,\cdot\cdot\cdot,e_d\}$ to $\{O,e_1^\ast,\cdot\cdot\cdot,e_d^\ast\}$ in $U_\alpha$, then  we have $e_i=\sum_{j=1}^da_{ij}e_j^\ast$. By (\ref{eq:mce2}) and (\ref{eq:ctm}), one easily obtains
\begin{equation}\label{eq:bcf}{\Omega_i^{\phantom{i}j}}^\ast=\sum_{k=1}^da_{ik}a_{jl}\Omega_k^{\phantom{k}l},\end{equation}
where ${\Omega_i^{\phantom{i}j}}^\ast$ are curvature forms defined by (\ref{eq:mce2}) with respect to the frame $\{O,e_1^\ast,\cdot\cdot\cdot e_d^\ast\}$. Substitute (\ref{eq:bcf}) into (\ref{eq:gbi}), $\Omega$ remains invariant. This shows that $\Omega_\alpha$ is intrinsic. By (\ref{eq:ltf}), $\Omega_\alpha=\Omega_\beta$ on $U_\alpha\bigcap U_\beta$, therefore the $\Omega_\alpha$s paste together to define a global form $\Omega$ on $M$.\qquad$\square$\\
\textsc{Remark}\\
By the above lemma, it follows immediate that
\begin{equation}\Omega=K(x)\omega^1\wedge\cdot\cdot\cdot\wedge\omega^d,K(x)\in C^\infty(M),\end{equation}
where $K(x)$ is a scalar invariant of the Riemannian manifold which is the generalization of the Gaussian curvature in the surface case. Therefore Lemma 3.2.1 should be viewed as a generalization of Gauss's \textit{Theorema Egregium} in high dimensions.\\
With these preliminaries we state our main theorem as follows:\\
\textbf{Theorem 3.2.2. (Gauss-Bonnet-Chern)} \textit{Let} $M$ \textit{be a closed oriented Riemannian manifold with an even dimension} $d$, $\Omega$ \textit{is defined as above}. \textit{Then the following formula holds}
\begin{equation}\label{eq:gbc44}\int_M\Omega=\chi(M).\end{equation}
\textsc{A Transgression Lemma}\\
When a closed form on $M$ is pulled back to a fiber bundle, it could happen that it becomes an exact form, such a process is called \textit{transgression}.\\
Let $\pi:SM\rightarrow M$ be the projective sphere bundle of $M$, then $\pi^\ast\Omega\in\Gamma(SM,\mathit{\Lambda}^dT^\ast SM)$.\\
\textbf{Lemma 3.2.3. (Transgression)} \textit{There exists a} $\Pi\in\Gamma(SM,\mathit{\Lambda}^{d-1}T^\ast SM)$ \textit{such that}
\begin{equation}\pi^\ast\Omega=d\Pi.\end{equation}
\textit{Proof}. Suppose $v\in\Gamma(U_\alpha,SM)$, then we have
\begin{equation}\label{eq:leng}\pi^\ast v=\pi^\ast(v^i e_i),\|v\|=\sqrt{\sum_{i=1}^dv_i^2}=1.\end{equation}
Apply exterior differentiation to $\pi^\ast v$, by (\ref{eq:ctm}) we get
\begin{equation}\label{eq:diff}d\pi^\ast v=\eta^i\pi^\ast e_i,\eta^i=d\pi^\ast v^i+\pi^\ast(v^j\omega_j^{\phantom{j}i}).\end{equation}
A dimensional consideration shows that the 1-forms $\{\eta^i\}$ are not independent. In fact, we can deduce from (\ref{eq:diff}) the important relation
\begin{equation}\label{eq:rela}\sum_{i=1}^d\pi^\ast v^i\eta^i=0.\end{equation}
Take the exterior differential of $\eta^i$ in (\ref{eq:diff}) and use (\ref{eq:mce2}) to get
\begin{eqnarray}
\label{eq:exp2}d\eta^i=\pi^\ast(dv^j\wedge\omega_j^{\phantom{j}i})+\pi^\ast(v^j d\omega_j^{\phantom{j}i})\qquad\qquad\qquad\nonumber\\
=\pi^\ast(dv^j\wedge\omega_j^{\phantom{j}i})+\pi^\ast\Big(v^j(\Omega_j^{\phantom{j}i}+\omega_j^{\phantom{j}k}\wedge\omega_k^{\phantom{k}i})\Big)\nonumber\\
=\pi^\ast(v^j\Omega_j^{\phantom{j}i})-\pi^\ast(\eta^i\wedge\omega_i^{\phantom{i}j}).\qquad\qquad\qquad\quad
\end{eqnarray}
Following Chern, we construct the following two sets of differential forms
\begin{eqnarray}
\Phi_k=\sum_{\sigma\in\Sigma_d}\textrm{sgn }\sigma\cdot \pi^\ast v^{\sigma(1)}\eta^{\sigma(2)}\wedge\cdot\cdot\cdot\wedge\eta^{\sigma(d-2k)}\wedge\nonumber\\
\pi^\ast\Omega_{\sigma(d-2k+1)}^{\phantom{\sigma(d-2k+1)}\sigma(d-2k+2)}\wedge\cdot\cdot\cdot\wedge\pi^\ast\Omega_{\sigma(d-1)}^{\phantom{\sigma(d-1)}\sigma(d)},
\end{eqnarray}
\begin{eqnarray}
\Psi_k=\sum_{\sigma\in\Sigma_d}\textrm{sgn }\sigma\cdot\pi^\ast\Omega_{\sigma(1)}^{\phantom{\sigma(1)}\sigma(2)}\wedge\eta^{\sigma(3)}\wedge\cdot\cdot\cdot\wedge\eta^{\sigma(d-2k)}\wedge\nonumber\\
\pi^\ast\Omega_{\sigma(d-2k+1)}^{\phantom{\sigma(d-2k+1)}\sigma(d-2k+2)}\wedge\cdot\cdot\cdot\wedge\pi^\ast\Omega_{\sigma(d-1)}^{\phantom{\sigma(d-1)}\sigma(d)},\qquad
\end{eqnarray}
where $0\leq k\leq \frac{d}{2}-1$. By a similar argument as in the proof of Lemma 3.2.1, it is easy to check that the differential forms $\Phi_k$ and $\Psi_k$ are intrinsic and therefore global by (\ref{eq:red}).\\
The main step for proving this lemma is to establish the following recurrent relation:
\begin{equation}\label{eq:rec}d\Phi_k=\Psi_{k-1}+\frac{d-2k+1}{2(k+1)}\Psi_k,\end{equation}
where $\Psi_{-1}$ is defined to be 0.\\
Now suppose (\ref{eq:rec}) has been proved, then we can solve $\Psi_k$ in terms of $d\Phi_0,\cdot\cdot\cdot,d\Phi_k$:
\begin{equation}\Psi_k=\sum_{m=0}^d(-1)^m\frac{2^{m+1}(k+1)!(d-2k+2m-2)!}{(k-m)!(d-2k-1)!}d\Phi_{k-m},\end{equation}
In particular, we get
\begin{equation}\pi^\ast\Omega=(-1)^{\frac{d}{2}-1}\frac{1}{2^d\pi^{\frac{d}{2}}(\frac{d}{2})!}\Psi_{\frac{d}{2}-1}=d\Pi,\end{equation}
where
\begin{equation}\Pi=\pi^{\frac{d}{2}}\sum_{m=0}^{\frac{d}{2}-1}(-1)^m\frac{(\frac{d}{2}-m-1)!}{(d-2m-1)!2^{\frac{d}{2}+m-1}}\Phi_m.\end{equation}
The proof is complete.\qquad$\square$\\
\textbf{Proof of (\ref{eq:rec})}\\
Computation gives
\begin{eqnarray}
d\Phi_k=\sum_{\sigma\in\Sigma_d}\textrm{sgn }\sigma\cdot\pi^\ast dv^{\sigma(1)}\wedge\eta^{\sigma(2)}\wedge\cdot\cdot\cdot\wedge\eta^{\sigma(d-2k)}\wedge\qquad\qquad\nonumber\\
\pi^\ast\Omega_{\sigma(d-2k+1)}^{\phantom{\sigma(d-2k+1)}\sigma(d-2k+2)}\wedge\cdot\cdot\cdot\wedge\pi^\ast\Omega_{\sigma(d-1)}^{\phantom{\sigma(d-1)}\sigma(d)}\qquad\qquad\qquad\nonumber\\
+(d-2k+1)\sum_{\sigma\in\Sigma_d}\textrm{sgn }\sigma\cdot\pi^\ast v^{\sigma(1)}d\eta^{\sigma(2)}\wedge\cdot\cdot\cdot\wedge\eta^{\sigma(d-2k)}\wedge\nonumber\\
\pi^\ast\Omega_{\sigma(d-2k+1)}^{\phantom{\sigma(d-2k+1)}\sigma(d-2k+2)}\wedge\cdot\cdot\cdot\wedge\pi^\ast\Omega_{\sigma(d-1)}^{\phantom{\sigma(d-1)}\sigma(d)}\qquad\qquad\qquad\nonumber\\
-k\sum_{\sigma\in\Sigma_d}\textrm{sgn }\sigma\cdot\pi^\ast v^{\sigma(1)}\eta^{\sigma(2)}\wedge\cdot\cdot\cdot\wedge\eta^{\sigma(d-2k)}\wedge\quad\qquad\qquad\nonumber\\
\pi^\ast d\Omega_{\sigma(d-2k+1)}^{\phantom{\sigma(d-2k+1)}\sigma(d-2k+2)}\wedge\cdot\cdot\cdot\wedge\pi^\ast\Omega_{\sigma(d-1)}^{\phantom{\sigma(d-1)}\sigma(d)}.\qquad\qquad\qquad
\end{eqnarray}
If we substitute $dv^i$, $d\eta^i$ and $d\Omega_i^{\phantom{i}j}$ by their expressions in (\ref{eq:diff}), (\ref{eq:exp2}) and (\ref{eq:bian}), the resulting expression for $d\Phi_k$ will then consist of two kinds of terms, those involving $\pi^\ast\omega_i^{\phantom{i}j}$ and those not. We collect the terms not involving $\pi^\ast\omega_i^{\phantom{i}j}$, which are
\begin{eqnarray}
\label{eq:Psi}\Psi_{k-1}+(d-2k+1)\sum_{\sigma\in\Sigma_d}\textrm{sgn }\cdot \pi^\ast(v^{\sigma(1)}v^j\Omega_j^{\phantom{j}\sigma(2)})\wedge\qquad\qquad\qquad\quad\nonumber\\
\eta^{\sigma(3)}\wedge\cdot\cdot\cdot\wedge\eta^{\sigma(d-2k)}\wedge\pi^\ast\Omega_{\sigma(d-2k+1)}^{\phantom{{\sigma(d-2k+1)}}{\sigma(d-2k+2)}}\wedge\cdot\cdot\cdot\wedge\pi^\ast\Omega_{\sigma(d-1)}^{\phantom{{\sigma(d-1)}}{\sigma(d)}}.
\end{eqnarray}
These expressions are clearly intrinsic and therefore global, actually it is easy to see that they are just $d\Phi_k$s.\\
In fact, fix an $O\in M$, we can choose a frame $\{O,e_1,\cdot\cdot\cdot,e_d\}$ at $O$ so that $\omega_i^{\phantom{i}j}=0$. This is equivalent to using the Riemannian normal coordinates in the tensor calculus context. Then it is obvious that (\ref{eq:Psi}) is equal to $d\Phi_k$ at $O$. It follows immediately they are identical because they are intrinsic.\\
In order to prove (\ref{eq:rec}), we use $\Psi_k$ to simplify the second term in (\ref{eq:Psi}). To do this, we make the following abbreviations
\begin{eqnarray}
A_k=\sum_{\sigma\in\Sigma_d}\textrm{sgn }\sigma\cdot\pi^\ast\Big((v^{\sigma(1)})^2\Omega_{\sigma(1)}^{\phantom{{\sigma(1)}}\sigma(2)}\Big)\wedge\eta^{\sigma(3)}\wedge\cdot\cdot\cdot\wedge\eta^{\sigma(d-2k)}\wedge\nonumber\\
\pi^\ast\Omega_{\sigma(d-2k+1)}^{\phantom{{\sigma(d-2k+1)}}{\sigma(d-2k+2)}}\wedge\cdot\cdot\cdot\wedge\pi^\ast\Omega_{\sigma(d-1)}^{\phantom{{\sigma(d-1)}}{\sigma(d)}},\qquad\qquad\qquad
\end{eqnarray}
\begin{eqnarray}
B_k=\sum_{\sigma\in\Sigma_d}\textrm{sgn }\sigma\cdot\pi^\ast(v^{\sigma(1)}v^{\sigma(3)}\Omega_{\sigma(3)}^{\phantom{{\sigma(3)}}\sigma(2)})\wedge\eta^{\sigma(3)}\wedge\cdot\cdot\cdot\wedge\eta^{\sigma(d-2k)}\wedge\nonumber\\
\pi^\ast\Omega_{\sigma(d-2k+1)}^{\phantom{{\sigma(d-2k+1)}}{\sigma(d-2k+2)}}\wedge\cdot\cdot\cdot\wedge\pi^\ast\Omega_{\sigma(d-1)}^{\phantom{{\sigma(d-1)}}{\sigma(d)}},\qquad\qquad\qquad\quad
\end{eqnarray}
\begin{eqnarray}
C_k=\sum_{\sigma\in\Sigma_d}\textrm{sgn }\sigma\cdot\pi^\ast\Big((v^{\sigma(3)})^2\Omega_{\sigma(1)}^{\phantom{{\sigma(1)}}\sigma(2)}\Big)\wedge\eta^{\sigma(3)}\wedge\cdot\cdot\cdot\wedge\eta^{\sigma(d-2k)}\wedge\nonumber\\
\pi^\ast\Omega_{\sigma(d-2k+1)}^{\phantom{{\sigma(d-2k+1)}}{\sigma(d-2k+2)}}\wedge\cdot\cdot\cdot\wedge\pi^\ast\Omega_{\sigma(d-1)}^{\phantom{{\sigma(d-1)}}{\sigma(d)}},\qquad\qquad\qquad
\end{eqnarray}
therefore we can write
\begin{equation}\label{eq:simp}d\Phi_k=\Psi_{k-1}+(d-2k-1)\Big(A_k+(d-2k-2)B_k\Big).\end{equation}
By (\ref{eq:leng}) we have
\begin{eqnarray}A_k=\sum_{\sigma\in\Sigma_d}\textrm{sgn }\sigma\cdot\pi^\ast\bigg(\Big(1-\sum_{i=2}^d(v^{\sigma(i)})^2\Big)\Omega_{\sigma(1)}^{\phantom{{\sigma(1)}}\sigma(2)}\bigg)\wedge\eta^{\sigma(3)}\wedge\cdot\cdot\cdot\wedge\eta^{\sigma(d-2k)}\nonumber\\
\wedge\pi^\ast\Omega_{\sigma(d-2k+1)}^{\phantom{{\sigma(d-2k+1)}}{\sigma(d-2k+2)}}\wedge\cdot\cdot\cdot\wedge\pi^\ast\Omega_{\sigma(d-1)}^{\phantom{{\sigma(d-1)}}{\sigma(d)}}\qquad\qquad\qquad\qquad\qquad\nonumber\\
=\Psi_k-A_k-(d-2k-2)C_k-2kA_k.\qquad\qquad\qquad\qquad\qquad\qquad\quad
\end{eqnarray}
This gives the relation
\begin{equation}\label{eq:re1}\Psi_k=2(k+1)A_k+(d-2k-2)C_k.\end{equation}
By (\ref{eq:rela}), we have
\begin{eqnarray}
\label{eq:re2}B_k=\sum_{\sigma\in\Sigma_d}\textrm{sgn }\sigma\cdot\pi^\ast(v^{\sigma(1)}\Omega_{\sigma(3)}^{\phantom{{\sigma(3)}}\sigma(2)})\wedge(-\sum_{i=1}^d\pi^\ast v^{\sigma(i)}\eta^{\sigma(i)})\wedge\qquad\qquad\qquad\nonumber\\
\eta^{\sigma(4)}\wedge\cdot\cdot\cdot\wedge\eta^{\sigma(d-2k)}\wedge\pi^\ast\Omega_{\sigma(d-2k+1)}^{\phantom{{\sigma(d-2k+1)}}{\sigma(d-2k+2)}}\wedge\cdot\cdot\cdot\wedge\pi^\ast\Omega_{\sigma(d-1)}^{\phantom{{\sigma(d-1)}}{\sigma(d)}}\nonumber\\
=C_k-(2k+1)B_k.\qquad\qquad\qquad\qquad\qquad\qquad\qquad\qquad\qquad\qquad\quad
\end{eqnarray}
Substitute (\ref{eq:re1}) and (\ref{eq:re2}) into (\ref{eq:simp}), we get (\ref{eq:rec}).\qquad$\square$\\
\textsc{Proof of the Gauss-Bonnet-Chern Theorem}\\
First note that $\mathcal{S}$ is a submanifold of $SM$, so we can use $\pi|_{\mathcal{S}}$ to pull the integral $\int_{M\setminus\bigcup_{x_i\in I}\overline{B(x_i)}}\Omega$ back to $\mathcal{S}$. By the transgression lemma and the Stokes' theorem, we have
\begin{eqnarray}
\label{eq:fin}\int_M\Omega=\int_{M\setminus\bigcup_{x_i\in I}\overline{B(x_i)}}\Omega+\int_{\bigcup_{x_i\in I}\overline{B(x_i)}}\Omega\qquad\nonumber\\
=\int_\mathcal{S}\pi|_\mathcal{S}^\ast\Omega+\int_{\bigcup_{x_i\in I}\overline{B(x_i)}}\Omega\qquad\qquad\quad\nonumber\\
=\int_\mathcal{S}d\Pi|_\mathcal{S}+\int_{\bigcup_{x_i\in I}\overline{B(x_i)}}\Omega\qquad\qquad\quad\nonumber\\
=\int_{\pi|_\mathcal{S}^{-1}\big(\bigcup_{x_i\in I}\partial\overline{B(x_i)}\big)}\Pi+\int_{\bigcup_{x_i\in I}\overline{B(x_i)}}\Omega.
\end{eqnarray}
As the radii of the balls $B(x_i)$ tend to 0, the first term on the right hand side of (\ref{eq:fin}) tends to $\chi(M)\int_{SM_x}\Pi$ by the Poincar\'{e}-Hopf index theorem and the second term vanishes, where $SM_x$ denotes the fiber of $SM$ at $x$.\\
By (\ref{eq:anti}), one easily calculates
\begin{equation}\int_{SM_x}\Pi=\frac{(\frac{d}{2})!\pi^\frac{d}{2}}{d!2^{\frac{d}{2}-1}}\int_{SM_x}\Phi_0=1.\end{equation} This completes the proof.\qquad$\square$\\
\textsc{Remarks}
\begin{itemize}
\item In Chern's original proof, he assumed that $\nabla^{TM}$ is the Levi-Civita connection. However, it is clear from the proof above that the Gauss-Bonnet-Chern theorem holds as long as $\nabla^{TM}$ is metric compatible.
\item Chern's method of transgression can be modified to generalize the Gauss-Bonnet-Chern theorem to odd dimensions, although the formula in this case is somewhat trivial and not of much interest. This is done by Chern himself in \cite{sc2}. In fact, define
\begin{equation}\widetilde{\Pi}=\frac{(-1)^d}{2^d\pi^{\frac{1}{2}(d-1)}}\sum_{k=0}^{[\frac{1}{2}(d-1)]}(-1)^k\frac{1}{k!\Gamma\Big(\frac{1}{2}(d-2k+1)\Big)}\widetilde{\Phi}_k,\end{equation}
where
\begin{eqnarray}
\widetilde{\Phi}_k=\sum_{\sigma\in\Sigma_d}\textrm{sgn }\sigma\cdot\pi^\ast\Omega_{\sigma(1)}^{\phantom{{\sigma(1)}}\sigma(2)}\wedge\cdot\cdot\cdot\wedge\pi^\ast\Omega_{\sigma(2k-1)}^{\phantom{{\sigma(2k-1)}}\sigma(2k)}\wedge\nonumber\\
\pi^\ast\omega_{\sigma(2k+1)}^{\phantom{{\sigma(2k+1)}}\sigma(d)}\wedge\cdot\cdot\cdot\wedge\pi^\ast\Omega_{\sigma(d-1)}^{\phantom{{\sigma(d-1)}}\sigma(d)}.\qquad\qquad\quad
\end{eqnarray}
Chern showed that $d\widetilde{\Pi}=\pi^\ast\widetilde{\Omega}$, where
\begin{equation}\widetilde{\Omega}=\left\{\begin{array}{ll}\Omega & d\textrm{ even},\\0&d\textrm{ odd}.\end{array}\right.\end{equation}
Therefore, the Gauss-Bonnet-Chern theorem for every closed oriented Riemannian manifold is given by
\begin{equation}\label{eq:gbco}\int_M\widetilde{\Omega}=\left\{\begin{array}{ll}\chi(M)&d\textrm{ even},\\0&d\textrm{ odd}.\end{array}\right.\end{equation}
Since $\chi(M)=0$ if $M$ is a closed manifold with an odd dimension, we actually have $\int_M\widetilde{\Omega}=\chi(M)$ for every closed oriented Riemannian manifold $M$.\\
With the forms $\widetilde{\Pi}$ and $\widetilde{\Omega}$, Chern obtained intrinsically in \cite{sc2} the Allendoerfer-Weil formula (c.f. (\ref{eq:awf})) for a differentiable Riemannian polyhedra $P$. Since $\partial P$ is a submanifold of $M$, let $\widetilde{\pi}:PN\partial P\bigcap SM\rightarrow \partial P$ be the bundle of inner unit normal vectors over $\partial P$, where $PN\partial P$ is the projectivized normal bundle, i.e., the quotient bundle of the normal bundle $N\partial P$ over $\partial P$ which treats $\Big\{(x,\lambda y):x\in \partial P,y\in N\partial P,\lambda\in\mathbb{R}\Big\}$ as a single point. Under these conventions the Allendoerfer-Weil formula is given by
\begin{equation}\int_P\widetilde{\Omega}=\int_{\partial P}(\widetilde{\pi}^{-1})^\ast\widetilde{\Pi}|_{PN\partial P\bigcap SM}-\chi'(M).\end{equation}
\end{itemize}
\subsection{Chern-Weil Theory}
After his excellent work on the Gauss-Bonnet theorem, Chern realized the possibility of expressing certain topological invariants by combinations of curvature forms. This leads to his fundamental paper \cite{sc6} introducing the \textit{Chern classes}
\begin{equation}\det\Big(I+\frac{\sqrt{-1}}{2\pi}\Theta(\nabla^E)\Big)=1+c_1\Big(\Theta(\nabla^E)\Big)+\cdot\cdot\cdot+c_r\Big(\Theta(\nabla^E)\Big),\end{equation}
where $I$ is the identity matrix, $E$ is a Hermitian vector bundle of rank $r$ over a differentiable manifold $M$, and $\Theta(\nabla^E)$ is the curvature matrix with respect to a connection $\nabla^E$ on $E$.\\
Chern showed in \cite{sc6} that the differential forms $c_i\Big(\Theta(\nabla^E)\Big),1\leq i\leq r$ are actually independent of $\nabla^E$ and determine cohomology classes $\bigg[c_i\Big(\Theta(\nabla^E)\Big)\bigg]$ in the de Rham cohomology groups. More precisely,
\begin{equation}\bigg[c_i\Big(\Theta(\nabla^E)\Big)\bigg]\in H^{2i}_{\textrm{dR}}(M,\mathbb{C}),1\leq i\leq r.\end{equation}
This example, together with the top form $\Omega$ defined by (\ref{eq:gbi}) in the proof of the Gauss-Bonnet-Chern theorem, give important examples of the Chern-Weil theory, which is a theory aimed at representing characteristic classes of some principal bundle over a manifold by geometrically significant objects. The above two examples due to Chern are unified in the Chern-Weil homomorphism, which is the main result of the Chern-Weil theory. To state this theorem, we need some preliminaries in multilinear algebra.\\
Let $\mathbb{K}=\mathbb{R}\textrm{ or }\mathbb{C}$. $G$ is a real or complex Lie group and $\mathfrak{g}$ its Lie algebra. We introduce a function $\nu:\mathfrak{g}\times\cdot\cdot\cdot\times\mathfrak{g}\rightarrow\mathbb{K}$ satisfying the following two properties:
\begin{itemize}
\item $\nu$ is $k$-linear and
\begin{equation}\nu(X_{\sigma(1)},\cdot\cdot\cdot,X_{\sigma(k)})=\nu(X_1,\cdot\cdot\cdot,X_k),\sigma\in\Sigma_k,\end{equation}
\item
\begin{equation}\nu\Big(\textrm{Ad}(g)X_1,\cdot\cdot\cdot,\textrm{Ad}(g)X_k\Big)=\nu(X_1,\cdot\cdot\cdot,X_k),g\in\mathfrak{g}.\end{equation}
\end{itemize}
It is a well-known fact that every such $\nu$ can be identified with a homogeneous polynomial of degree $k$ defined by $\nu(X)=\nu(X,\cdot\cdot\cdot,X)$, $X\in\mathfrak{g}$. We will call $\nu(X)$ an \textit{invariant polynomial}. It is clear that all the invariant polynomials form a ring, which we will denote by $I(G)$. Note that the action of a $\nu\in I(G)$ can be extended to $\mathfrak{g}$-valued differential forms in the following manner:
\begin{equation}\nu(\zeta^i\otimes Y_i)=\zeta^i\otimes\nu(Y_i),\zeta^i\in\mathit{\Lambda}^\ast T^\ast M,Y_i\in\mathfrak{g}.\end{equation}
\textbf{Theorem 3.3.1. (Chern-Weil Homomorphsim)} \textit{For a} $d$-\textit{dimensional differentiable manifold} $M$, \textit{let} $\nabla^P$ \textit{be a Cartan connection on a principal} $G$-\textit{bundle} $\pi:P\rightarrow M$, \textit{then} $\nabla^P$ \textit{is given by a set of} $\mathfrak{g}$-\textit{valued} 1-\textit{forms} $\Big\{\theta_\alpha(\nabla^P)\Big\}$, \textit{let} $\Big\{\Theta_\alpha(\nabla^P)\Big\}$ \textit{be the corresponding curvature} 2-\textit{forms}. $\nu\in I(G)$. \textit{Then} $\nu\Big(\Theta_\alpha(\nabla^P)\Big)$ \textit{is intrinsic and thus global by} (\ref{eq:ltf}), \textit{we have}
\begin{itemize}
\item $\nu\Big(\pi^\ast\Theta_\alpha(\nabla^P)\Big)$ \textit{is a coboundary in} $P$, \textit{i}.\textit{e}., \textit{there is a} $\xi\in\Gamma(P,\mathit{\Lambda}^{d-1}T^\ast M)$ \textit{such that} $\nu\Big(\pi^\ast\Theta_\alpha(\nabla^P)\Big)=d\xi$ on $P$. \textit{This is a generalization of the transgression lemma} 3.2.3 \textit{proved above}.
\item $\nu\Big(\Theta_\alpha(\nabla^P)\Big)$ \textit{is closed in} $M$, \textit{and} $\bigg[\nu\Big(\Theta_\alpha(\nabla^P)\Big)\bigg]\in H_{\textrm{dR}}^{2k}(M,\mathbb{K})$ \textit{is independent of} $\nabla^P$, \textit{therefore this defines a ring homomorphsim}
\begin{equation}\mathcal{W}:I(G)\rightarrow H_{\textrm{dR}}^\ast(M,\mathbb{K})\end{equation}
\textit{by taking} $\nu$ to $\bigg[\nu\Big(\Theta_\alpha(\nabla^P)\Big)\bigg]$. $\mathcal{W}$ \textit{is called the Chern-Weil homomorphism}.
\end{itemize}
For a proof of this theorem and a neat and profound introduction to the theory of characteristic classes, we refer the reader to the appendix: \textit{Geometry of Characteristic Classes} of \cite{sc4}. See also \cite{jms} for a topological approach to characteristic classes.\\
We now use this theorem to give another interpretation of the Gauss-Bonnet-Chern theorem. This interpretation leads naturally to the topological proofs in the next section.\\
\textbf{Definition 3.3.2.} Let $A=(a_{ij})$ be a $d\times d$ skew-symmetric matrix, $d$ even, then the \textit{Pfaffian} of $A$ is defined by
\begin{equation}\textrm{Pf}(A)=\frac{1}{2^\frac{d}{2}(\frac{d}{2})!}\sum_{\sigma\in\Sigma_d}\Big(\textrm{sgn }\sigma\cdot\prod_{i=1}^\frac{d}{2}a_{\sigma(2i-1)\sigma(2i)}\Big).\end{equation}
One is able to verify the following important relation:
\begin{equation}\label{eq:pfbp}\textrm{Pf}(A)^2=\det(A).\end{equation}
Using this notation, we can write the Gauss-Bonnet integrand as
\begin{equation}\Omega=(\frac{-1}{2\pi})^\frac{d}{2}\textrm{Pf}\Big(\Theta(\nabla^{TM})\Big),\end{equation}
where $\nabla^{TM}$ denotes a metric connection on $TM$. Since $\textrm{Pf}\in I\Big(\textrm{O}(d,\mathbb{R})\Big)$, by the Chern-Weil homomorphsim, $\bigg[\textrm{Pf}\Big(\Theta(\nabla^{TM})\Big)\bigg]\in H_{\textrm{dR}}^d(M,\mathbb{R})$ is independent of the metric connection $\nabla^{TM}$ we choose, so we can simply write $\Omega=(\frac{-1}{2\pi})^\frac{d}{2}\textrm{Pf}(\Theta^{TM})$.\\
The $d$-form $(\frac{-1}{2\pi})^\frac{d}{2}\textrm{Pf}(\Theta^{TM})$ is usually called the \textit{geometric Euler class}. Then the Gauss-Bonnet-Chern theorem asserts that \textit{the integration of the geometric Euler class on a closed oriented Riemannian manifold gives the Euler characteristic}.
\section{Topological Proofs}

\subsection{Thom Class and Euler Class}
To make our exposition self-contained, we shall provide here a brief introduction to the topological backgrounds which will be assumed in our following proofs. Our main object in this subsection is to prove the following simple but important result:\\
\textbf{Proposition 4.1.1.} \textit{The pullback of the Thom class of an even dimensional oriented vector bundle via the zero section to the base manifold is the Euler class}.\\
We present here a new approach to the Euler class inspired by \cite{bt} and \cite{mt}, as we shall see, this approach is very convenient for our purpose and doesn't seem to be too formal to hide the geometric motivation. Compare our exposition here with \cite{agp}, \cite{bt}, \cite{ah}, \cite{mt} and \cite{jms}.\\
\textsc{Thom Class}\\
As we know, there are two basic topological invariants for a manifold $M$, namely the de Rham cohomology $H^\ast_{\textrm{dR}}(M,\mathbb{R})$ and the compactly supported cohomology $H^\ast_{\textrm{c}}(M,\mathbb{R})$. For a vector bundle $\pi:E\rightarrow M$, there is a third one, namely the cohomology with compact support in the vertical direction, we use $\Omega_{\textrm{cv}}^k(E)$ to denote the corresponding complex.\\
\textbf{Definition 4.1.2.} A smooth $k$-form $\omega\in\Omega_{\textrm{cv}}^k(E)$ if and only if for every compact set $K\subset M$, $\pi^{-1}(K)\bigcap\textrm{Supp }\omega$ is compact. The cohomology of the complex $\Omega_{\textrm{cv}}^\ast(E)$, denoted by $H_{\textrm{cv}}^\ast(E,\mathbb{R})$, is called the \textit{cohomology of} $E$ \textit{with compact support in the vertical direction}.\\
Just as the Poincar\'{e} duality which gives the relation between $H^\ast_{\textrm{dR}}(M,\mathbb{R})$ and $H^\ast_{\textrm{c}}(M,\mathbb{R})$, the Thom isomorphism relates $H^\ast_{\textrm{cv}}(E,\mathbb{R})$ with $H^\ast_{\textrm{dR}}(M,\mathbb{R})$.\\
\textbf{Theorem 4.1.3. (Thom Isomorphism)} \textit{Let} $E$ \textit{be an oriented vector bundle over a manifold} $M$ \textit{with rank} $n$, \textit{then}
\begin{equation}H_{\textrm{cv}}^\ast(E,\mathbb{R})\cong H_{\textrm{dR}}^{\ast-n}(M,\mathbb{R}).\end{equation}
The proof is based on the generalized Mayer-Vietoris argument and the Poincar\'{e} lemma for compactly supported cohomology, details can be found in \cite{agp}, \cite{bt} and \cite{mt}.\\
By this theorem $1\in H_{\textrm{dR}}^0(M,\mathbb{R})$ determines a cohomology class $t(E)\in H_{\textrm{cv}}^n(E,\mathbb{R})$, we call $t(E)$ the \textit{Thom class} of the oriented vector bundle $E$. There is a simple characterization of the Thom class:\\
\textbf{Proposition 4.1.4.} \textit{The Thom class} $t(E)$ \textit{on a rank} $n$ \textit{oriented vector bundle} $E$ \textit{can be uniquely characterized as the cohomology class in} $H_{\textrm{cv}}^n(E,\mathbb{R})$ \textit{which restricts to the generator of} $H_{\textrm{c}}^n(E_x,\mathbb{R})$ \textit{for every} $x\in M$.\\
\textit{Proof}. By Definition 4.1.4, $\textrm{Supp }t(E)|_{E_x}$ is compact. Since the pushforward $\pi_\ast$ is given by integration along the fiber and $\pi_\ast t(E)=1$, we see that $t(E)|_{E_x}$ is a generator of $H_{\textrm{c}}^n(E_x,\mathbb{R})$.\\
Conversely if $t'(E)\in H_{\textrm{cv}}^n(E,\mathbb{R})$ restricts to a generator of $H_{\textrm{c}}^n(E_x,\mathbb{R})$ on each fiber, then one verifies directly the so-called \textit{projection formula}:
\begin{equation}\pi_\ast\Big((\pi^\ast\omega)\wedge t'(E)\Big)=\omega\wedge\pi_\ast t'(E)=\omega.\end{equation}
Since it is another simple verification that the map $\mathfrak{T}:\omega\rightarrow(\pi^\ast\omega)\wedge t'(E)$ defines the Thom isomorphism, it follows that $t'(E)=\mathfrak{T}(1)$ defines the Thom class.\qquad$\square$\\
\textbf{Corollary 4.1.5.} \textit{Suppose} $E_1$ \textit{is an orientable vector bundle over} $M$ \textit{with rank} $n$ \textit{and} $f:N\rightarrow M$ \textit{is} $C^\infty$. \textit{Let} $E_2$ \textit{be a vector bundle over} $N$ \textit{such that the following diagram}
\[\begin{CD}E_2@>\tilde{f}>>E_1\\
@VV\tilde{\pi} V @VV\pi V\\
N@>f>>M
\end{CD}\]
\textit{is commutative, then} $E_2$ \textit{must be orientable}. \textit{We have} $\tilde{f}^\ast t(E_1)=t(E_2)$.\\
\textit{Proof}. It is clear that $\tilde{f}^\ast t(E_1)\in H_{\textrm{cv}}^n(E_2)$. By the commutativity of the diagram,
\begin{equation}\tilde{\pi}_\ast\tilde{f}^\ast t(E_1)=f^\ast\pi_\ast t(E_1)=f^\ast1=1.\end{equation}
By Proposition 4.1.4, the proof is complete.\qquad$\square$\\
\textsc{Euler Class}\\
We begin with the simplest case when $E$ is an oriented Riemannian plane bundle and $\{U_\alpha\}$ is a coordinate open cover of $M$. Then the local coordinates in $E|_{U_\alpha}$ can be written as $\pi^\ast x_1,\cdot\cdot\cdot,\pi^\ast x_n,r_\alpha,\vartheta_\alpha$, where $x_1,\cdot\cdot\cdot,x_n$ is the local coordinates on $U_\alpha$ and every fiber of $E$ is endowed with the polar coordinates. Since the structure group of $E$ can be reduced to SO(2), we may assume that on $U_\alpha\bigcap U_\beta$, $r_\alpha=r_\beta$ while $\vartheta_\alpha$ and $\vartheta_\beta$ differ by a rotation. Define $\varphi_{\alpha\beta}$, up to a constant multiple of $2\pi$, as the angle of rotation in the counterclockwise direction from the $\alpha$-coordinate system to the $\beta$-coordinate system:
\begin{equation}\label{eq:the}\vartheta_\beta=\vartheta_\alpha+\pi^\ast\varphi_{\alpha\beta},\varphi_{\alpha\beta}:U_\alpha\bigcap U_\beta\rightarrow\mathbb{R}.\end{equation}
By the definition of $\varphi_{\alpha\beta}$, we have the simple relation
\begin{equation}\varphi_{\alpha\beta}+\varphi_{\beta\gamma}+\varphi_{\alpha\gamma}\in2\pi\mathbb{Z},\end{equation}
it follows immediately that $d\varphi_{\alpha\beta}$ satisfies the cocycle condition
\begin{equation}\label{eq;coc}d\varphi_{\alpha\beta}+d\varphi_{\beta\gamma}+d\varphi_{\alpha\gamma}=0.\end{equation}
Let $\{\rho_\gamma\}$ be a partition of unity subordinate to $\{U_\gamma\}$, and define $\xi_\alpha=\frac{1}{2\pi}\sum_\gamma\rho_\gamma d\varphi_{\alpha\gamma}$. By (\ref{eq;coc}) we have
\begin{equation}\label{eq:pie}\frac{1}{2\pi}d\varphi_{\alpha\beta}=\xi_\beta-\xi_\alpha.\end{equation}
(\ref{eq:pie}) shows that $d\xi_\alpha$ agrees with $d\xi_\beta$ on $U_\alpha\bigcap U_\beta$, hence $d\xi_\alpha$ piece together to give a global 2-form $e^E(x)$ on $M$. Since $e^E(x)$ is closed, it determines a cohomology class $e(E)$ in $H_{\textrm{dR}}^2(M,\mathbb{R})$, which is by definition the \textit{Euler class} of the oriented plane bundle $E$. It is obvious that $e(E)$ does not depend on the choice of $\xi_\alpha$.\\
A primary reason for considering the plane bundle case is that the Euler class can be expressed explicitly by transition funcitons in this case. This makes it easier for us to obtain the basic properties of the Euler class in the plane bundle case. In fact, let $g_{\alpha\beta}:U_\alpha\bigcap U_\beta\rightarrow\textrm{SO}(2)$ denote the transition functions of $E$, then by the identification
\begin{equation}\left(\begin{array}{ll}\cos\vartheta & -\sin\vartheta\\ \sin\vartheta & \quad\cos\vartheta\end{array}\right)\rightarrow e^{i\vartheta},\end{equation}
$g_{\alpha\beta}$ may be thought of as complex-values functions. In this context, we can write
\begin{equation}-\pi^\ast\varphi_{\alpha\beta}=\vartheta_\alpha-\vartheta_\beta=-i\pi^\ast\log g_{\alpha\beta}.\end{equation}
By the injectivity of $\pi^\ast$ we get
\begin{equation}\label{eq:varp}\varphi_{\alpha\beta}=-i\log g_{\alpha\beta}.\end{equation}
By (\ref{eq:varp}) and the definition of $\xi_\alpha$,
\begin{equation}\xi_\alpha=-\frac{1}{2\pi i}\sum_\gamma\rho_\gamma d\log g_{\gamma\alpha}.\end{equation}
Apply exterior differentiation to $\xi_\alpha$, we finally get
\begin{equation}\label{eq:tran}e^E(x)=-\frac{1}{2\pi i}\sum_\gamma d(\rho_\gamma d\log g_{\gamma\alpha})\textrm{ on }U_\alpha.\end{equation}
Identify $M$ with the image of the zero section of $E$. By (\ref{eq:the}) and (\ref{eq:pie}), we have on $(E\setminus M)|_{U_\alpha\bigcap U_\beta}$
\begin{equation}\frac{d\vartheta_\alpha}{2\pi}-\pi^\ast\xi_\alpha=\frac{d\vartheta_\beta}{2\pi}-\pi^\ast\xi_\beta.\end{equation}
These forms then piece together to give a global 1-form $\psi^E(y)$ on $E\setminus M$. One verifies easily that
\begin{equation}d\psi^E(y)=-\pi^\ast e^E(x).\end{equation}
Let $\rho(r)\in C_0^\infty(\mathbb{R}_+\bigcup\{0\})$ be a function which is identically $-1$ near $r=0$. Define
\begin{equation}\mathsf{t}(E)=\bigg[d\Big(\rho(r)\cdot\psi^E(y)\Big)\bigg].\end{equation}
An easy verification shows that $\mathsf{t}(E)$ satisfies the following properties
\begin{itemize}
\item $\mathsf{t}(E)\subset\Omega_{\textrm{cv}}^2(E)$;
\item $\pi_\ast\iota_x^\ast\mathsf{t}(E)=1$, where $\iota_x:E_x\rightarrow E$ is the inclusion;
\item $\mathsf{t}(E)\in H_{\textrm{cv}}^n(E,\mathbb{R})$ is independent of the choice of $\rho(r)$.
\end{itemize}
By Proposition 4.1.4, we conclude that $\mathsf{t}(E)=t(E)$. Let $s_0:M\rightarrow E$ be the zero section of $E$, we have
\begin{equation}\label{eq:spe}s_0^\ast\mathsf{t}(E)=\bigg[d\Big(\rho(0)\Big)\cdot s_0^\ast\psi^E(y)-\rho(0)s_0^\ast\pi^\ast e^E(x)\bigg]=e(E).\end{equation}
Therefore we have proved Proposition 4.1.1 in the plane bundle case.\\
We now use the following result to extend the definition of the Euler class to an arbitrary even-dimensional oriented vector bundle $E$.\\
\textbf{Theorem 4.1.6. (Splitting Principle for Real Vector Bundles)} \textit{Let} $E$ \textit{be an oriented vector bundle over a manifold} $M$ \textit{with rank} $2n$. \textit{Then there is a manifold} $N$ \textit{and a map} $\mathcal{J}:N\rightarrow M$ \textit{such that}
\begin{itemize}
\item \textit{the homomorphism} $\mathcal{J}^\ast:H_{\textrm{dR}}^\ast(M,\mathbb{R})\rightarrow H_{\textrm{dR}}^\ast(N,\mathbb{R})$ is \textit{injective};
\item $\mathcal{J}^\ast E$ \textit{is a direct sum of orientable plane bundles}: $\mathcal{J}^\ast E=\bigoplus_{i=1}^n\mathcal{P}_i$.
\end{itemize}
A proof of this theorem can be found in \cite{mt} and \cite{ps}, see also \cite{bg} for a general version of the real splitting principle.\\
\textbf{Definition 4.1.7.} Let $E$ be a $2n$-dimensional oriented vector bundle over a manifold $M$, the \textit{Euler class} $e(E)$ is defined by the following formula
\begin{equation}\label{eq:eul}\mathcal{J}^\ast e(E)=e(\mathcal{P}_1)\smile\cdot\cdot\cdot\smile e(\mathcal{P}_n),\end{equation}
where $\smile$ denotes the cup product.\\
This Euler class shall sometimes be referred to as the \textit{topological Euler class} to distinguish it from the geometric Euler class defined in $\S3.3$ since it is defined in a topological way.\\
\textsc{Remark}\\
By now it is not clear whether the Euler class in the above definition is well-defined, especially the existence of $e(E)$ in (\ref{eq:eul}) needs to be proved. However, we will show in the following proposition that $s_0^\ast t(E)$ satisfies (\ref{eq:eul}). Combined with (\ref{eq:spe}), it is easy to see that this will establish Proposition 4.1.1 and the existence of $e(E)$ simultaneously. Note that the injectivity of $\mathcal{J}^\ast$ assures the uniqueness of $e(E)$ as long as it exists.\\
\textbf{Proposition 4.1.8.}
\begin{equation}\mathcal{J}^\ast(s_0^\ast t)(E)=(s_0^\ast t)(\mathcal{J}^\ast E).\end{equation}
\textit{Proof}. We have the following commutative diagram
\[\begin{CD}\mathcal{J}^\ast E @>\tilde{\mathcal{J}}>> E\\
@VV\tilde{\pi} V @VV\pi V\\
N @>\mathcal{J}>> M.
\end{CD}\]
Let $s_0$ and $\tilde{s}_0$ be the zero sections of $E$ and $\mathcal{J}^\ast E$ respectively, then we have $\tilde{\mathcal{J}}\circ \tilde{s}_0=s_0\circ\mathcal{J}$. By Corollary 4.1.5, the proof is complete.\qquad$\square$\\
\textsc{Euler Number and Euler Characteristic}\\
Let $E$ be a rank $2n$ oriented vector bundle over $M$. Since $H_{\textrm{dR}}^{2n}(M,\mathbb{R})=\mathbb{R}$, the Euler class may be identified with the number $\int_Me(E)$, which is by definition the \textit{Euler number} of $E$. The following result shows that the concept of the Euler number is a generalization of the Euler characteristic.\\
\textbf{Proposition 4.1.9.} \textit{Suppose} $M$ \textit{is closed and oriented and $\dim_{\mathbb{R}}M=d$, $d$ even. Then}
\begin{equation}\label{eq:enec}\int_Me(TM)=\chi(M).\end{equation}
\textit{Proof}. Let $\{\omega_i\}$ be a basis of the vector space $H_{\textrm{dR}}^\ast(M,\mathbb{R})$, $\{\tau_j\}$ be the dual basis under the Poincar\'{e} duality. Let $\rho_1$ and $\rho_2$ be two projections from $M\times M$ onto $M$. By the K\"{u}nneth formula, $H_{\textrm{dR}}^\ast(M\times M,\mathbb{R})=H_{\textrm{dR}}^\ast(M,\mathbb{R})\otimes H_{\textrm{dR}}^\ast(M,\mathbb{R})$, it follows that $\{\rho_1^\ast\omega_i\wedge\rho_2^\ast\tau_j\}$ forms a basis of $H_{\textrm{dR}}^\ast(M\times M,\mathbb{R})$. Suppose $\eta_\blacktriangle$ is the Poincar\'{e} dual of the diagonal $\blacktriangle\in M\times M$, then we can write
\begin{equation}\eta_\blacktriangle=c^{ij}\rho_1^\ast\omega_i\wedge\rho_2^\ast\tau_j,c^{ij}\in\mathbb{R}.\end{equation}
\textbf{Lemma 4.1.10.}
\begin{equation}\eta_\blacktriangle=(-1)^{\textrm{deg }\omega_i}\rho_1^\ast\omega_i\wedge\rho_2^\ast\tau_i.\end{equation}
\textbf{Proof of Lemma 4.1.10}\\
We compute $\int_\blacktriangle\rho_1^\ast\tau_k\wedge\rho_2^\ast\omega_l$ in two ways. On one hand, use the map $\iota: M\rightarrow\blacktriangle\subset M\times M$ to pull this integral back to $M$, we get
\begin{equation}\label{eq:fir}\int_\blacktriangle\rho_1^\ast\tau_k\wedge\rho_2^\ast\omega_l=\int_M\iota^\ast\rho_1^\ast\tau_k\wedge\iota^\ast\rho_2^\ast\omega_l=\int_M\tau_k\wedge\omega_l=(-1)^{\textrm{deg }\tau_k\cdot\textrm{deg }\omega_l}\delta_{kl}.\end{equation}
On the other hand, by the definition of the Poincar\'{e} dual, we compute
\begin{eqnarray}
\label{eq:seco}\int_\blacktriangle\rho_1^\ast\tau_k\wedge\rho_2^\ast\omega_l=\int_{M\times M}\rho_1^\ast\tau_k\wedge\rho_2^\ast\omega_l\wedge\eta_\blacktriangle\qquad\qquad\qquad\qquad\qquad\qquad\quad\nonumber\\
=c^{ij}\int_{M\times M}\rho_1^\ast\tau_k\wedge\rho_2^\ast\omega_l\wedge\rho_1^\ast\omega_i\wedge\rho_2^\ast\tau_j\qquad\qquad\qquad\qquad\nonumber\\
=(-1)^{\textrm{deg }\omega_i(\textrm{deg }\tau_k+\textrm{deg }\omega_l)}c^{ij}\int_{M\times M}\rho_1^\ast(\omega_i\wedge\tau_k)\rho_2^\ast(\omega_l\wedge\tau_j)\nonumber\\
=(-1)^{\textrm{deg }\omega_k(\textrm{deg }\tau_k+\textrm{deg }\omega_l)}c^{kl}.\qquad\qquad\qquad\qquad\qquad\qquad
\end{eqnarray}
Compare (\ref{eq:fir}) with (\ref{eq:seco}), the proof is complete.\qquad$\square$\\
A straightforward verification shows the following diagram commutes:
\[\begin{CD}0 @>>> T\blacktriangle @>>> T\Big((M\times M)|_\blacktriangle\Big) @>>> N\blacktriangle @>>> 0\\
@. @| @|\\
0 @>>> TM @>>> TM\oplus TM @>>> TM @>>> 0,
\end{CD}\]
it then follows that
\begin{equation}\label{eq:isom}N\blacktriangle\cong TM\cong T\blacktriangle.\end{equation}
It is a well-known fact that the Poincar\'{e} dual of a closed oriented submanifold $N$ is represented by the Thom class of a tubular neighborhood (which is isomorphic to the normal bundle) of $N$ (see \cite{bt} for detail). Use this we calculate
\begin{equation}\label{eq:intp}\int_\blacktriangle\eta_\blacktriangle=\int_\blacktriangle t(N\blacktriangle)=\int_\blacktriangle e(N\blacktriangle)=\int_\blacktriangle e(T\blacktriangle)=\int_Me(TM),\end{equation}
where Proposition 4.1.1 and (\ref{eq:isom}) have been used during the calculation.\\
On the other hand, by Lemma 4.1.10,
\begin{eqnarray}\label{eq:pd}\int_\blacktriangle\eta_\blacktriangle=(-1)^{\textrm{deg }\omega_i}\int_\blacktriangle\rho_1^\ast\omega_i\wedge\rho_2^\ast\tau_i\qquad\nonumber\\
=(-1)^{\textrm{deg }\omega_i}\int_M\iota^\ast\rho_1^\ast\omega_i\wedge\iota^\ast\rho_2^\ast\tau_i\nonumber\\
=(-1)^{\textrm{deg }\omega_i}\int_M\omega_i\wedge\tau_i\qquad\quad\nonumber\\
=\sum_{j=0}^d(-1)^j\dim H_{\textrm{dR}}^j(M,\mathbb{R})\quad\nonumber\\
=\chi(M).\qquad\qquad\qquad\qquad\quad
\end{eqnarray}
Combine (\ref{eq:intp}) and (\ref{eq:pd}), we get (\ref{eq:enec}).\qquad$\square$
\subsection{A Generalized Gauss-Bonnet-Chern Theorem}
It is clear from our interpretation of the Gauss-Bonnet-Chern theorem at the end of $\S3.3$ and Proposition 4.1.9 that to prove the Gauss-Bonnet-Chern theorem, we only need to show the geometric Euler class and the topological Euler class of $TM$ coincide, i.e., $\Big[(\frac{-1}{2\pi})^\frac{d}{2}\textrm{Pf}(\Theta^{TM})\Big]=e(TM)$.\\
It is an observation made by Denis Bell in \cite{db} that it is essentially not more difficult to prove the general fact:
\begin{equation}\label{eq:ggbc}\Big[(\frac{-1}{2\pi})^\frac{d}{2}\textrm{Pf}(\Theta^E)\Big]=e(E),\end{equation}
where $E$ is an oriented Riemannian vector bundle of an even rank $d$ and $\Theta^E$ is its curvature matrix induced by a metric connection $\nabla^E$. By the Chern-Weil homomorphism, $\Big[\textrm{Pf}(\Theta^E)\Big]$ is independent of the $\nabla^E$ we choose. Evidently, (\ref{eq:ggbc}) generalizes the Gauss-Bonnet-Chern theorem.\\
\textbf{Lemma 4.2.1.} \textit{Suppose} $E$ \textit{is a Riemannian vector bundle over a manifold} $M$, \textit{then} $(\frac{-1}{2\pi})^\frac{d}{2}\textrm{Pf}(\Theta^E)$ \textit{satisfies the following two properties}:
\begin{itemize}
\item \textit{When} $E=E_1\oplus\cdot\cdot\cdot\oplus E_n$, \textit{then}
\begin{equation}(\frac{-1}{2\pi})^\frac{d}{2}\textrm{Pf}(\Theta^E)=(\frac{-1}{2\pi})^\frac{d_1}{2}\textrm{Pf}(\Theta^{E_1})\wedge\cdot\cdot\cdot\wedge(\frac{-1}{2\pi})^\frac{d_n}{2}\textrm{Pf}(\Theta^{E_n}),\end{equation}
\textit{where} $d_i=\textrm{rank }E_i$ \textit{and} $d=\sum_{i=1}^nd_i$.
\item \textit{Let} $f:N\rightarrow M$ \textit{be a} $C^\infty$ \textit{map between two manifolds}, \textit{then}
\begin{equation}\label{eq:pbpf}(\frac{-1}{2\pi})^\frac{d}{2}\textrm{Pf}(\Theta^{f^\ast E})=f^\ast\Big((\frac{-1}{2\pi})^\frac{d}{2}\textrm{Pf}(\Theta^E)\Big).\end{equation}
\end{itemize}
\textit{Proof}. Without loss of generality we can assume $n=2$. Suppose $\theta(\nabla^{E_1})$ and $\theta(\nabla^{E_2})$ are connection matrices for $E_1$ and $E_2$ respectively, then it is clear that
\begin{equation}\theta(\nabla^E)=\left(\begin{array}{ll}\theta(\nabla^{E_1}) & 0 \\ 0 & \theta(\nabla^{E_2})\end{array}\right)\end{equation}
defines a connection matrix for $E_1\oplus E_2$. Moreover, if $\theta(\nabla^{E_1})$ and $\theta(\nabla^{E_2})$ are skew-symmetric, so is $\theta(\nabla^E)$.\\
Therefore the corresponding curvature matrix for $E_1\oplus E_2$ is given by
\begin{eqnarray}
\Theta^E=\left(\begin{array}{ll}d\theta(\nabla^{E_1})-\theta(\nabla^{E_1})\wedge\theta(\nabla^{E_1}) & \quad\qquad\qquad0 \\ \quad\qquad\qquad0 & d\theta(\nabla^{E_2})-\theta(\nabla^{E_2})\wedge\theta(\nabla^{E_2})\end{array}\right)\nonumber\\
=\left(\begin{array}{ll}\Theta^{E_1} & 0 \\ 0 & \Theta^{E_2}\end{array}\right).\qquad\qquad\qquad\qquad\qquad\qquad\qquad\qquad\qquad\qquad
\end{eqnarray}
By Definition 3.3.2, this gives
\begin{equation}(\frac{-1}{2\pi})^\frac{d}{2}\textrm{Pf}(\Theta^E)=(\frac{-1}{2\pi})^\frac{d_1}{2}\textrm{Pf}(\Theta^{E_1})\wedge(\frac{-1}{2\pi})^\frac{d_2}{2}\textrm{Pf}(\Theta^{E_2}).\end{equation}
Since it is obvious that the pullback of a skew-symmetric connection matrix is still a skew-symmetric connection matrix on the pullback vector bundle, (\ref{eq:pbpf}) holds trivially.\qquad$\square$\\
By this lemma, in order to prove (\ref{eq:ggbc}), we only need to establish (\ref{eq:ggbc}) in the plane bundle case. The general case then follows from Definition 4.1.7 and the above lemma. Therefore, the generalized Gauss-Bonnet-Chern theorem reduces to the following proposition:\\
\textbf{Proposition 4.2.2.} \textit{Suppose} $\pi:\mathcal{P}\rightarrow M$ \textit{is an oriented Riemannian plane bundle}, \textit{then} $\Big[-\frac{1}{2\pi}\textrm{Pf}(\Theta^\mathcal{P})\Big]=e(\mathcal{P})$.\\
\textit{Proof}. First note that in this case, $\textrm{Pf}(\Theta^\mathcal{P})$ is given by $d{(\theta_\alpha)}_1^{\phantom{1}2}$ on $U_\alpha$, where ${(\theta_\alpha)}_1^{\phantom{1}2}$ is the connection 1-form for an arbitrary metric connection on $\mathcal{P}$. This fact is familiar to us when the connection $\nabla^\mathcal{P}$ is the Levi-Civita connection on $TM$, where $M$ is a Riemannian surface (cf. \cite{ccl}). In the general case it can be easily verified by a direct computation.\\
The proof will be finished as long as we can show (\ref{eq:tran}) and $-\frac{1}{2\pi}d{(\theta_\alpha)}_1^{\phantom{1}2}$ defines the same cohomology class in $H_{\textrm{dR}}^2(M,\mathbb{R})$ for some metric connection $\nabla^{\mathcal{P}}$.\\
Since $\mathcal{P}$ is oriented, we have $g_{\alpha\beta}\in\textrm{SO}(2)$. By identifying $\textrm{SO}(2)$ with $\mathbb{S}^1$, we can write $g_{\alpha\beta}=e^{i\varphi_{\alpha\beta}}$ for some $\varphi_{\alpha\beta}\in(0,2\pi]$. By (\ref{eq:varp}), it is evident that this $\varphi_{\alpha\beta}$ is identical to that defined in (\ref{eq:the}). We then substitute $g_{\gamma\alpha}$ by $e^{i\varphi_{\gamma\alpha}}$ in (\ref{eq:tran}), this gives
\begin{eqnarray}
\label{eq:locf}e(\mathcal{P})|_{U_\alpha}=\Big[-\frac{1}{2\pi i}\sum_\gamma d(\rho_\gamma d\log g_{\gamma\alpha})\Big]\nonumber\\
=\Big[-\frac{1}{2\pi i}\sum_\gamma d(\rho_\gamma g_{\gamma\alpha}^{-1}dg_{\gamma\alpha})\Big]\nonumber\\
=\Big[-\frac{1}{2\pi}\sum_\gamma d(\rho_\gamma d\varphi_{\gamma\alpha})\Big].\quad
\end{eqnarray}
We now associate to each fiber of $\mathcal{P}$ an orthogonal frame $\{O,\mathbf{r}_\alpha,\mathbf{r}_\alpha^\bot\}$ such that $\{\pi^\ast O,\pi^\ast\mathbf{r}_\alpha,\pi^\ast\mathbf{r}_\alpha^\bot\}$ forms the Cartesian coordinate system of the fiber and $\|\pi^\ast\mathbf{r}_\alpha\|=\|\pi^\ast\mathbf{r}_\alpha^\bot\|=1$ (identify the tautological sections with vectors in $\mathcal{P}$). This is equivalent to considering the special orthogonal frame bundle $\mathcal{F}$ of $\mathcal{P}$. By (\ref{eq:ctm}), the connection 1-forms of $\mathcal{F}$ are given by
\begin{equation}\label{eq:con2}d\mathbf{r}_\alpha={(\theta_\alpha)}_1^{\phantom{1}2}\mathbf{r}_\alpha^\bot.\end{equation}
On the other hand, by (\ref{eq:the}) we have
\begin{equation}d\pi^\ast\mathbf{r}_\alpha=\Big(\pi^\ast d\varphi_{\gamma\alpha}+\pi^\ast{(\theta_\gamma)}_1^{\phantom{1}2}\Big)\pi^\ast\mathbf{r}_\alpha^\bot.
\end{equation}
By (\ref{eq:con2}) and the injectivity of $\pi^\ast$, we have on $U_\alpha\bigcap U_\gamma$:
\begin{equation}\label{eq:crit}{(\theta_\alpha)}_1^{\phantom{1}2}=d\varphi_{\gamma\alpha}+{(\theta_\gamma)}_1^{\phantom{1}2}.\end{equation}
Substitute (\ref{eq:crit}) into (\ref{eq:locf}), we get
\begin{eqnarray}e(\mathcal{P})|_{U_\alpha}=\Bigg[-\frac{1}{2\pi}\sum_\gamma d\bigg(\rho_\gamma\Big({(\theta_\alpha)}_1^{\phantom{1}2}-{(\theta_\gamma)}_1^{\phantom{1}2}\Big)\bigg)\Bigg]\qquad\nonumber\\
=\Big[-\frac{1}{2\pi}d\sum_\gamma\rho_\gamma{(\theta_\alpha)}_1^{\phantom{1}2}+\frac{1}{2\pi}d\sum_\gamma\rho_\gamma{(\theta_\gamma)}_1^{\phantom{1}2}\Big]\nonumber\\
=\Big[-\frac{1}{2\pi}d{(\theta_\alpha)}_1^{\phantom{1}2}+\frac{1}{2\pi}d\sum_\gamma\rho_\gamma{(\theta_\gamma)}_1^{\phantom{1}2}\Big]\qquad\quad\nonumber\\
=\Big[-\frac{1}{2\pi}d{(\theta_\alpha)}_1^{\phantom{1}2}\Big].\qquad\qquad\qquad\qquad\qquad\quad
\end{eqnarray}
The last step in the above computation is because $\frac{1}{2\pi}d\sum_\gamma\rho_\gamma{(\theta_\gamma)}_1^{\phantom{1}2}$ defines an exact 2-form on $M$, and therefore a coboundary in the de Rham complex. This proves that $e(\mathcal{P})=\Big[-\frac{1}{2\pi}\textrm{Pf}(\Theta^\mathcal{P})\Big]\in H_{\textrm{dR}}^2(M,\mathbb{R})$.\qquad$\square$\\
In particular, we have the following version of the Gauss-Bonnet-Chern theorem:\\
\textbf{Theorem 4.2.3. (Generalized Gauss-Bonnet-Chern Theorem)} \textit{Suppose} $E$ \textit{is a rank} $d$ \textit{oriented Riemannian vector bundle over a} $d$-\textit{dimensional closed oriented manifold} $M$, $d$ \textit{even}. \textit{Then we have the following formula}:
\begin{equation}(\frac{-1}{2\pi})^\frac{d}{2}\int_M\textrm{Pf}(\Theta^E)=\int_Me(E).\end{equation}
\textsc{Remarks}
\begin{itemize}
\item Such a theorem was first appeared in \cite{ms} as a consequence of the representation of the Pontrjagin classes by curvature forms. However, since they didn't make use of the real splitting principle, they restrict themselves to the case when $E$ is a direct sum of plane bundles. When $E=TM$, a topological proof based on the relation between the Euler class and the first Chern class can be found in \cite{mt}. The above version of the Gauss-Bonnet-Chern theorem was first proved in \cite{db}, the idea is to construct the Thom class of $\mathcal{P}$ geometrically. Our new proof above is a direct approach, which I believe is the simplest proof for this theorem.
\item Note that our approach to proving the Gauss-Bonnet-Chern theorem actually shows the equivalence between the classical Gauss-Bonnet theorem on surfaces (\ref{eq:ori1}) and the Gauss-Bonnet-Chern theorem for all even dimensions (\ref{eq:gbc44}). However, this does not diminish the importance of Chern's work \cite{sc1} since we may never come to realize this equivalence without his contributions to the whole subject.
\end{itemize}
Theorem 4.2.3 enables us to generalize the classical Poincar\'{e}-Hopf index theorem as follows:\\
\textbf{Theorem 4.2.4. (Generalized Poincar\'{e}-Hopf Index Theorem)} \textit{Let} $E$ \textit{be an oriented rank} $d$ \textit{vector bundle over a closed oriented} $d$-\textit{dimensional manifold} $M$, $d$ \textit{even}. $S$ \textit{is the sphere bundle of} $E$. $\mathfrak{S}\in\Gamma(M\setminus I,S)$, \textit{where} $I$ \textit{is chosen to be isolated}. \textit{Denote by} $\textrm{locdeg}_\mathfrak{S}(x)$ \textit{the local degree of} $\mathfrak{S}$ \textit{at} $x\in I$, \textit{then the following formula holds}:
\begin{equation}\label{eq:gph}\sum_{x_i\in I}\textrm{locdeg}_\mathfrak{S}(x_i)=\int_Me(E).\end{equation}
\textit{Proof}. It is an easy observation that every step before (\ref{eq:fin}) in Chern's proof in $\S3.2$ can actually be carried out on an arbitrary Riemannian vector bundle $E$ of rank $d$. Endow $E$ with a Riemannian structure, by applying Chern's approach to $E$ instead of $TM$, we get
\begin{equation}(\frac{-1}{2\pi})^\frac{d}{2}\int_M\textrm{Pf}(\Theta^E)=\int_{\pi|^{-1}_\mathfrak{S}\big(\bigcup_{x_i\in I}\partial\overline{B(x_i)}\big)}\Pi_S+(\frac{-1}{2\pi})^\frac{d}{2}\int_{\bigcup_{x_i\in I}\overline{B(x_i)}}\textrm{Pf}(\Theta^E),\end{equation}
where $\Pi_S$ is the corresponding $(d-1)$-form constructed on $\pi:S\rightarrow M$ such that $(\frac{-1}{2\pi})^\frac{d}{2}\pi^\ast\textrm{Pf}(\Theta^E)=d\Pi_S$.\\
Now let the radii of $B(x_i)$ tend to 0, apply Theorem 4.2.3 to deduce (\ref{eq:gph}).\qquad$\square$\\
\textsc{Remark}\\
For a different proof of this theorem, see \cite{bt}. Such a theorem actually holds for every oriented sphere bundle with fiber $\mathbb{S}^{d-1}$.
\subsection{Berezin Integral and A Model of Construction}
The reason for passing to the Thom class to carry out the required geometric construction lies in Proposition 4.1.4, i.e., since the Thom class is easily characterized by its properties restricted on each fiber, we may hope that our construction would be completed by some suitable modifications on a relatively trivial construction for a single fiber. This is actually what we plan to do in the following text.\\
As a preparation for the Thom form proof in the next subsection, we introduce here the concept of a Berezin integral, which will be used as a technical tool in our geometric construction of Mathai-Quillen's Thom form. As an application of such a tool, we will use the Berezin integral to construct the Thom form for the vector bundle over a point, this is a model of the general case and provides enough motivation for the general construction. We follow from \cite{bgv}, \cite{wz1} and \cite{wz}.\\
\textbf{Definition 4.3.1.} Let $V$ be a real vector space. A nonzero linear map $B:\mathit{\Lambda}^\ast V\rightarrow\mathbb{R}$ which vanishes on $\mathit{\Lambda}^kV$ for $k<\dim_{\mathbb{R}}V$ is called a \textit{Berezin integral}, by $\mathit{\Lambda}^\ast V$ we mean the exterior algebra of $V$.\\
To understand this definition, we restrict ourselves to the simplest case when $V$ is an oriented Euclidean space with basis $\{e_1,\cdot\cdot\cdot,e_n\}$. Then there is a canonical Berezin integral given by projecting $\omega\in\mathit{\Lambda}^\ast V$ onto $e_1\wedge\cdot\cdot\cdot\wedge e_n$, i.e.,
\begin{equation}B(e^I)=\left\{\begin{array}{ll}1 & |I|=n,\\0 & \textrm{otherwise},\end{array}\right.\end{equation}
where $I$ denotes the multi-index $I=\{i_k|1\leq k\leq n,i_1<\cdot\cdot\cdot<i_n\}$. It follows from the above definition that $B(\omega)\neq0$ if and only if the component of degree $n$ of $\omega\in\mathit{\Lambda}^\ast V$ is not 0.\\
Since the identity map $\textrm{Id}:V\rightarrow V$ can be identified with an element of $\Gamma(V,\mathit{\Lambda}^0 V^\ast\otimes V)$, we can take its exterior differential $dx\in\Gamma(V,\mathit{\Lambda}^1V^\ast\otimes V)$. Then the exponential $e^{-idx}$ lies in $\Gamma(V,\mathit{\Lambda}^\ast V^\ast\otimes\mathit{\Lambda}^\ast V)$. We extend the Berezin integral to a map $B:\Gamma(V,\mathit{\Lambda}^\ast V^\ast\otimes\mathit{\Lambda}^\ast V)\rightarrow\Gamma(V,\mathit{\Lambda}^\ast V^\ast)$ by
\begin{equation}B(\omega\otimes\xi)=\omega B(\xi),\omega\in\Gamma(V,\mathit{\Lambda}^\ast V^\ast),\xi\in\mathit{\Lambda}^\ast V.\end{equation}
Consider $V$ as a vector bundle over a point, by Proposition 4.1.4, a Thom form $t^V(x)\in t(V)$ is a compactly supported $n$-form on $V$ with $\int_Vt^V(x)=1$. Set
\begin{equation}u^V(x)=(2\pi)^{-\frac{n}{2}}e^{-\frac{1}{2}\|x\|^2}dx^1\wedge\cdot\cdot\cdot\wedge dx^n.\end{equation}
Then it is easy to verify that
\begin{equation}\int_Vu^V(x)=1,\end{equation}
which means $u^V(x)$ is a Thom form on $V$ except that $u^V(x)$ is of exponential decay instead of having a compact support, we will show how to remedy this flaw in the next subsection. We now interpret $u^V(x)$ using the Berezin integral.\\
\textbf{Proposition 4.3.2. }
\begin{equation}\label{eq:toy}u^V(x)=\varepsilon(n)(2\pi)^{-\frac{n}{2}}B(e^{-\frac{\|x\|^2}{2}-idx}),\end{equation}
\textit{where}
\begin{equation}\varepsilon(n)=\left\{\begin{array}{ll}1 & n\textrm{ even,}\\i & n\textrm{ odd.}\end{array}\right.\end{equation}
\textit{Proof}. Let $\{e_k\}$ be the dual basis of $dx^k$, we have
\begin{eqnarray}
B(e^{-idx})=B\Big(\prod_{k=1}^n(1-idx^k\otimes e_k)\Big)\qquad\qquad\qquad\quad\nonumber\\
=(-i)^nB\Big((dx^1\otimes e_1)\wedge\cdot\cdot\cdot\wedge(dx^n\otimes e_n)\Big)\nonumber\\
=(-i)^n(-1)^{\frac{n(n-1)}{2}}dx^1\wedge\cdot\cdot\cdot\wedge dx^n.\qquad\quad
\end{eqnarray}
This proves the proposition.\qquad$\square$\\
\textsc{Remark}\\
One can imagine how to modify (\ref{eq:toy}) when the vector space $V$ is replaced by an oriented Riemannian vector bundle $p:\mathcal{V}\rightarrow M$. Since a metric connection $\nabla^\mathcal{V}$ generalizes the exterior differential, we may replace $dx$ by $\nabla^\mathcal{V}x\in\Gamma(\mathcal{V},\mathit{\Lambda}^1T^\ast\mathcal{V}\otimes p^\ast\mathcal{V})$, where $p^\ast\mathcal{V}$ is a vector bundle over $\mathcal{V}$ and $x$ is treated as the tautological section of $p^\ast\mathcal{V}$. However, since $\nabla^\mathcal{V}$ may not be flat, it seems that some further modifications on $u^V$ will be needed.
\subsection{Mathai-Quillen's Thom Form}
Suppose $\mathcal{E}$ is an oriented Riemannian vector bundle of rank $n$ over the base manifold $M$, then we can extend the Berezin intergal once again to a map
\begin{equation}B:\Gamma(M,\mathit{\Lambda}^\ast T^\ast M\otimes\mathit{\Lambda}^\ast\mathcal{E})\rightarrow\Gamma(M,\mathit{\Lambda}^\ast T^\ast M)\end{equation}
by an obvious fiberwise extension.\\
For convenience, we shall use here the language of covariant derivatives to specify the notion of a connection on vector bundles. This enables us to work with global geometric objects during our geometric construction of a global Thom form.\\
Note that the covariant derivative $\nabla^\mathcal{E}:\Gamma(M,\mathcal{E})\rightarrow\Gamma(M,T^\ast M\otimes\mathcal{E})$ induces in a natural way a covariant derivative on $\mathit{\Lambda}^\ast\mathcal{E}$, which we still denote by $\nabla^\mathcal{E}$.\\
\textbf{Proposition 4.4.1.} \textit{Suppose} $\nabla^\mathcal{E}$ \textit{is metric compatible}, \textit{then for every} $\eta\in\Gamma(M,\mathit{\Lambda}^\ast T^\ast M\otimes\mathit{\Lambda}^\ast\mathcal{E})$, \textit{we have}
\begin{equation}dB(\eta)=B(\nabla^\mathcal{E}\eta).\end{equation}
\textit{Proof}. By definition of the Berezin integral, $B$ can be regarded as a section $B\in\Gamma(M,\mathit{\Lambda}^n\mathcal{E}^\ast)$ paired to the sections of $\mathit{\Lambda}^n\mathcal{E}$. Suppose $\{e^i\}$ is a local frame of $\mathcal{E}^\ast$, then $B$ is locally given by $e^1\wedge\cdot\cdot\cdot\wedge e^n$. The proposition follows from the Leibniz rule since $\nabla^\mathcal{E}B=0$.\qquad$\square$\\
We apply this proposition with $M=\mathcal{V}$ and $\mathcal{E}=p^\ast\mathcal{V}$, and with connection $\nabla^{p^\ast\mathcal{V}}=p^\ast\nabla^\mathcal{V}$. $\Gamma\Big(\mathcal{V},\mathit{\Lambda}^\ast T^\ast\mathcal{V}\otimes\mathit{\Lambda}^\ast(p^\ast\mathcal{V})\Big)$ is a bigraded algebra and admits the following decomposition:
\begin{equation}\Gamma\Big(\mathcal{V},\mathit{\Lambda}^\ast T^\ast\mathcal{V}\otimes\mathit{\Lambda}^\ast(p^\ast\mathcal{V})\Big)=\bigoplus_{i,j=1}^n\mathcal{A}^{i,j},\end{equation}
where \begin{equation}\mathcal{A}^{i,j}=\Gamma\Big(\mathcal{V},\mathit{\Lambda}^iT^\ast\mathcal{V}\otimes\mathit{\Lambda}^j(p^\ast\mathcal{V})\Big).\end{equation}
The covariant derivative defines a map $\nabla^{p^\ast\mathcal{V}}:\mathcal{A}^{i,j}\rightarrow\mathcal{A}^{i+1,j}$.\\
For $\omega\in\Gamma(\mathcal{V},\mathit{\Lambda}^\bullet T^\ast\mathcal{V})$, the contraction by $s\in\Gamma(\mathcal{V},p^\ast\mathcal{V})=\mathcal{A}^{0,1}$ is defined by
\begin{equation}a(s)\Big(\omega\otimes(s_1\wedge\cdot\cdot\cdot\wedge s_j)\Big)=\sum_{k=1}^j(-1)^{\textrm{deg }\omega+k-1}\langle s,s_k\rangle\omega\otimes(s_1\wedge\cdot\cdot\cdot\wedge\widehat{s_k}\wedge\cdot\cdot\cdot\wedge s_j),\end{equation}
where $s_k\in\Gamma(\mathcal{V},p^\ast\mathcal{V})$. This defines a map $a(s):\mathcal{A}^{i,j}\rightarrow\mathcal{A}^{i,j-1}$. In the following text, we shall identify the Lie algebra bundle $\mathfrak{so}(\mathcal{V})$ with $\mathit{\Lambda}^2\mathcal{V}$, this identification can be made explicitly by the following map
\begin{equation}\label{eq:ide}A\in\mathfrak{so}(\mathcal{V})\rightarrow\sum_{i<j}\langle Ae_i,e_j\rangle e_i\wedge e_j.\end{equation}
Since for any $s\in\mathcal{A}^{0,1}$ and $\eta\in\bigoplus_{i,j=1}^n\mathcal{A}^{i,j}$, $B\Big(a(s)\eta\Big)=0$, we have by Proposition 4.4.1 the formula
\begin{equation}\label{eq:bfb}dB(\eta)=B\bigg(\Big(\nabla^{p^\ast\mathcal{V}}-ia(s)\Big)\eta\bigg).\end{equation}
Let us explain some objects needed in our geometric construction of Mathai-Quillen's Thom form:
\begin{itemize}
\item the tautological section $x\in\mathcal{A}^{0,1}$;
\item the norm $\|x\|^2\in\mathcal{A}^{0,0}$ and the covariant derivative $\nabla^{p^\ast\mathcal{V}}x\in\mathcal{A}^{1,1}$;
\item identifying the curavture form $\Theta^\mathcal{V}=(\nabla^\mathcal{V})^2\in\Gamma\Big(M,\mathit{\Lambda}^2T^\ast M\otimes\mathfrak{so}(\mathcal{V})\Big)$ with an element of $\Gamma(M,\mathit{\Lambda}^2T^\ast M\otimes\mathit{\Lambda}^2\mathcal{V})$ and pulling it back to $\mathcal{V}$ using the projection $p$, we obatain an element $\Theta^{p^\ast\mathcal{V}}\in\mathcal{A}^{2,2}$.
\end{itemize}
\textbf{Lemma 4.4.2.}
\begin{itemize}
\item \textit{Let} $\mathcal{Q}=\frac{\|x\|^2}{2}+i\nabla x+\Theta^{p^\ast\mathcal{V}}\in\bigoplus_{i,j=1}^n\mathcal{A}^{i,j}$. \textit{Then}
\begin{equation}\Big(\nabla^{p^\ast\mathcal{V}}-ia(x)\Big)\mathcal{Q}=0.\end{equation}
\item \textit{If} $\rho\in C^\infty(\mathbb{R})$, \textit{define} $\rho(\mathcal{Q})\in\bigoplus_{i,j=1}^n\mathcal{A}^{i,j}$ \textit{by the formula}
\begin{equation}\rho(\mathcal{Q})=\sum_{k=0}^n\frac{\rho^{(k)}\Big(\frac{\|x\|^2}{2}\Big)}{k!}(i\nabla^{p^\ast\mathcal{V}}+\Theta^{p^\ast\mathcal{V}})^k,\end{equation}
\textit{then}
\begin{equation}\label{eq:lfc}\Big(\nabla^{p^\ast\mathcal{V}}-ia(x)\Big)\rho(\mathcal{Q})=0.\end{equation}
\end{itemize}
\textit{Proof}. By definition of the covariant derivative, we have $\nabla^{p^\ast\mathcal{V}}\Big(\|x\|^2\Big)=-2a(x)\nabla^{p^\ast\mathcal{V}}x$. Applying the covariant derivative twice gives the curvature $\nabla^{p^\ast\mathcal{V}}(\nabla^{p^\ast\mathcal{V}}x)=a(x)\Theta^{p^\ast\mathcal{V}}$. While the Bianchi identity says that $\nabla^{p^\ast\mathcal{V}}\Theta^{p^\ast\mathcal{V}}=0$. Combining these facts we have
\begin{eqnarray}
\Big(\nabla^{p^\ast\mathcal{V}}-ia(x)\Big)\mathcal{Q}=-a(x)\nabla^{p^\ast\mathcal{V}}x+ia(x)\Theta^{p^\ast\mathcal{V}}+a(x)\nabla^{p^\ast\mathcal{V}}x-ia(x)\Theta^{p^\ast\mathcal{V}}\nonumber\\
=0,\qquad\qquad\qquad\qquad\qquad\qquad\qquad\qquad\qquad\qquad\qquad
\end{eqnarray}
where we have used the obvious fact $a(x)\|x\|^2=0$.\\
Since $\nabla^{p^\ast\mathcal{V}}-ia(x)$ is a derivation of the algebra $\bigoplus_{i,j=1}^n\mathcal{A}^{i,j}$, (\ref{eq:lfc}) follows.\qquad$\square$\\
Since $\mathcal{Q}\in\bigoplus_{0\leq 0\leq 2}\mathcal{A}^{k,k}$, it follows that $\rho(\mathcal{Q})\in\bigoplus_{0\leq 0\leq n}\mathcal{A}^{k,k}$, therefore $B\Big(\rho(\mathcal{Q})\Big)\in\Gamma(\mathcal{V},\mathit{\Lambda}^nT^\ast\mathcal{V})$. By (\ref{eq:bfb}) and the above lemma, we get
\begin{equation}dB\Big(\rho(\mathcal{Q})\Big)=B\bigg(\Big(\nabla^{p^\ast\mathcal{V}}-ia(x)\Big)\rho(\mathcal{Q})\bigg)=0,\end{equation}
this shows that $B\Big(\rho(\mathcal{Q})\Big)$ is a closed $n$-form on $\mathcal{V}$.\\
We are now at a point to define \textit{Mathai-Quillen's Thom form} on $\mathcal{V}$, this is done by choosing $\rho$ to be $(2\pi)^{-\frac{n}{2}}\varepsilon(n)e^{-x}$ in $B\Big(\rho(\mathcal{Q})\Big)$:
\begin{equation}\label{eq:mqt}u^\mathcal{V}=(2\pi)^{-\frac{n}{2}}\varepsilon(n)B(e^{-\mathcal{Q}})=(2\pi)^{-\frac{n}{2}}\varepsilon(n)B\Big(e^{-\big(\frac{\|x\|^2}{2}+i\nabla^{p^\ast\mathcal{V}}x+\Theta^{p^\ast\mathcal{V}}\big)}\Big).\end{equation}
The following proposition shows that the form $u^\mathcal{V}$ constructed above has the same nature as its model $u^V$ constructed in $\S4.3$.\\
\textbf{Proposition 4.4.3.}
\begin{equation}\int_{\mathcal{V}/M}u^\mathcal{V}=1.\end{equation}
\textit{Proof}. To calculate the above integal, it suffices to consider the case in which $M$ is a point, but this is what Proposition 4.3.1 said.\qquad$\square$\\
\textsc{Remark}\\
Mathai-Quillen's Thom form constructed above has rapid decay at infinity instead of having compact supports on each fiber. However, using diffeomorphism from the interior of the unit ball bundle $\mathcal{BV}$ of $\mathcal{V}$ to $\mathcal{V}$ given by
\begin{equation}\upsilon:y\rightarrow\frac{y}{\Big(1-\|y\|^2\Big)^\frac{1}{2}},\end{equation}
we can pull $u^\mathcal{V}$ back to obtain a Thom form $\upsilon^\ast u^\mathcal{V}$ with support in $\mathcal{BV}$.\\
\textsc{Apply to the Gauss-Bonnet-Chern Theorem}\\
We are now prepared to prove the Gauss-Bonnet-Chern theorem using our geometric construction (\ref{eq:mqt}). Assume $\mathcal{V}$ is an oriented rank $n$ Riemannian vector bundle over an oriented $n$-dimensional closed manifold $M$, $s\in\Gamma(M,\mathcal{V})$. We have
\begin{equation}s^\ast\upsilon^\ast u^\mathcal{V}=(2\pi)^{-\frac{n}{2}}\varepsilon(n)B\Big(e^{-\big(\frac{\|s\|^2}{2}+i\nabla^{\mathcal{V}}s+\Theta^\mathcal{V}\big)}\Big).\end{equation}
In particular, when $s$ is the zero section and $n$ is even, by Proposition 4.1.1, we get a geometric expression of the Euler class:
\begin{equation}\label{eq:gcel}e(\mathcal{V})=(2\pi)^{-\frac{n}{2}}B(e^{\Theta^\mathcal{V}}).\end{equation}
\textsc{Remark}\\
It follws from (\ref{eq:gcel}) that the diffeomorphsim $\upsilon$ is superfluous, i.e., $s^\ast u^\mathcal{V}=s^\ast\upsilon^\ast u^\mathcal{V}$. This shows that the pullback by the zero section of a larger class than the Thom class is the Euler class, or equivalently, the converse of Proposition 4.1.1 is not true. Therefore Mathai-Quillen's Thom form should be viewed as an extension of the classical Thom class with the pullback property preserved.\\
\textbf{Lemma 4.4.4.} \textit{For every} $A\in\mathit{\Lambda}^2V$, \textit{where} $V$ \textit{is a vector space}, \textit{we deonte by} $\exp A$ \textit{its exponential in} $\mathit{\Lambda}^\ast V$. \textit{Then we have the following identity}:
\begin{equation}\textrm{Pf}(A)=B(\exp A).\end{equation}
\textit{Proof}. Under the identification (\ref{eq:ide}), we have
\begin{equation}B(\exp A)=B\Big(\exp\sum_{i<j}\langle Ae_i,e_j\rangle e_i\wedge e_j\Big).\end{equation}
Using the definition of a Pfaffian, a direct computation yields
\begin{equation}\textrm{Pf}(A)=B\Big(\exp\sum_{i<j}\langle Ae_i,e_j\rangle e_i\wedge e_j\Big).\end{equation}
This finishes the proof.\qquad$\square$\\
By this lemma, (\ref{eq:gcel}) becomes $e(\mathcal{V})=(2\pi)^{-\frac{n}{2}}\textrm{Pf}(\Theta^\mathcal{V})$. Integrate on $M$ gives:
\begin{equation}(2\pi)^{-\frac{n}{2}}\int_M\textrm{Pf}(\Theta^\mathcal{V})=\int_Me(\mathcal{V}).\end{equation}
This is just the Gauss-Bonnet-Chern theorem we proved in $\S4.2$.\\
\textsc{Remark}\\
This proof of the Gauss-Bonnet-Chern theorem first appeared in \cite{mq} as a consequence of their elegant construction (\ref{eq:mqt}). Their proof made use of the transgression of $u^\mathcal{V}$ and the Poincar\'{e}-Hopf index theorem, see also \cite{wz1} and \cite{wz} for expositions of their proof. Our exposition here is a simplification of their proof.
\section{A Heat Equation Approach}

\subsection{Prologue}
Intuitively, the heat flow on a Riemannian manifold $M$ is locally controlled by the curvature tensor. That's why the short time behavior of the heat operator can be used to study the geometric structures of $M$. A treatise on these applications is contained in \cite{bgm}. On the other hand, when $t\rightarrow\infty$, the heat flows over the whole manifold and therefore provides global information of the space, i.e., the topological invariants of $M$. The most important result on the long time behavior of the heat operator is the well-known theorem of Hodge concerning harmonic forms, see \cite{pg}, \cite{jj}, and \cite{sr}. Therefore it is natural to expect that by forming certain time-independent combinations of the local components of the heat kernel, one is able to link the geometric information with the topological invariants, thus formulating a proof of the Gauss-Bonnet-Chern theorem.\\
In this section, we shall give a heat equation proof of the Gauss-Bonnet-Chern theorem. Although the proof here may not be as elegant as the intrinsic argument which we discussed in $\S3$, the heat equation method is of great power and it can be developed to prove the Atiyah-Singer index theorem (see \cite{bgv} and \cite{pg} for expositions). Since this method has the feature of linking the local information with the global information, it is nothing strange that the proof here does not depend on the Poincar\'{e}-Hopf index theorem.\\
Our exposition here follows from \cite{pg}, \cite{vp} and \cite{sr}.
\subsection{Construction of the Heat Kernel}
We introduce the basic notions and construct the heat kernel for differential forms in this subsection. Throughout this section, $(M,g)$ will be a closed oriented Riemannian manifold with an even dimension $d$.\\
\textbf{Definition 5.2.1.} Let $V$ be a vector bundle over $M$ and $C^\infty(M,V)$ be the space of its smooth global sections. Let $P:C^\infty(M,V)\rightarrow C^\infty(M,V)$ be a self-adjoint elliptic pseudodifferential operator with a positive definite leading symbol. The \textit{heat equation} is the following partial differential equation:
\begin{equation}\label{eq:heat}(\frac{\partial}{\partial t}+P)f(t,x)=0,t\geq 0,f(0,x)=f(x).\end{equation}
\textsc{Remark}\\
Since we want to prove here the Gauss-Bonnet-Chern Theorem 3.3.2, it is therefore appropriate to restrict ourselves to the case of the heat operator for the Laplacian on differential forms $\frac{\partial}{\partial t}+\Delta$. Here $\Delta$ denotes the geometers' Laplacian given by $d\delta+\delta d$, $\delta$ being the adjoint of $d$.\\
\textbf{Definition 5.2.2.} A \textit{heat kernel} is a double form $e(t,x,y)\in C^\infty(\mathbb{R}_+\times M\times M,\mathit{\Lambda}^\bullet T^\ast M\otimes\mathit{\Lambda}^\bullet T^\ast M)$ satisfying the following conditions:
\begin{itemize}
\item \begin{equation}\label{eq:hk1}(\frac{\partial}{\partial t}+\Delta_y)e(t,x,y)=0,\end{equation}
where $\Delta_y$ means the Laplacian acts on $y$;
\item \begin{equation}\label{eq:hk2}\lim_{t\rightarrow 0}\int_Me(t,x,y)\wedge\star\omega(y)=\omega(x),\forall \omega\in L^2(M,\mathit{\Lambda}^\bullet T^\ast M),\end{equation}
where $\star$ denotes the Hodge star operator and $L^2(M,\mathit{\Lambda}^\bullet T^\ast M)$ is the space of all $L^2$-sections of the vector bundle $\mathit{\Lambda}^\bullet T^\ast M$ with respect to the metric $\langle\omega,\eta\rangle=\star(\omega\wedge\star\eta)$, where $\eta\in L^2(M,\mathit{\Lambda}^\bullet T^\ast M)$.
\end{itemize}
\textsc{Remark}\\
It is clear from the above definition that \begin{equation}\omega(t,x)=\int_Me(t,x,y)\wedge\star\omega(y)\end{equation}
is the solution of (\ref{eq:heat}) with $P$ replaced by $\Delta_x$. This shows that to construct the heat kernel is equivalent to solving the heat equation. However, since a variable $y$ has been added to ``decompose" (\ref{eq:heat}) into the above two equations (\ref{eq:hk1}) and (\ref{eq:hk2}), it is more convenient for us to study the heat kernel since it enables us to treat the equation and the initial condition separately on $y$ and $x$.\\
An intermediate step for constructing the heat kernel this is to construct the parametrix, which is characterized in the following definition.\\
\textbf{Definition 5.2.3.} A \textit{parametrix} for the heat operator $\frac{\partial}{\partial t}+\Delta_y$ is a double form $H(t,x,y)\in C^\infty(\mathbb{R}_+\times M\times M,\mathit{\Lambda}^\bullet T^\ast M\otimes\mathit{\Lambda}^\bullet T^\ast M)$ which satisfies the following conditions:
\begin{itemize}
\item $(\frac{\partial}{\partial t}+\Delta_y)H\in C(\mathbb{R}_+\bigcup\{0\}\times M\times M,\mathit{\Lambda}^\bullet T^\ast M\otimes\mathit{\Lambda}^\bullet T^\ast M)$;
\item $\lim_{t\rightarrow 0}\int_M H(t,x,y)\wedge\star\omega(y)=\omega(x),\omega\in C^{\infty}(U,\mathit{\Lambda}^\bullet T^\ast M)$.
\end{itemize}
For differential $p$-forms, let $N>\frac{d}{2}$ and define
\begin{equation}\label{eq:hnp}H_N^p(t,x,y)=(4\pi t)^{-\frac{d}{2}}e^{-\frac{r^2(x,y)}{4t}}\sum_{i=0}^{N}t^iu^{i,p}(x,y),i,p\in\mathbb{Z},(x,y)\in U\times U,\end{equation}
where $r(x,y)$ is the geodesic distance between $x$ and $y$, $U\times U$ is a sufficiently small open neighborhood of the diagonal $\blacktriangle\subset M\times M$ and $u^{i,p}(x,y)\in C^{\infty}(U\times U,\mathit{\Lambda}^pT^\ast M\otimes\mathit{\Lambda}^pT^\ast M)$ are to be determined.\\
Note first that every $u^{i,p}(x,y)$ can be naturally identified with an element in $\textrm{Hom}_{\mathbb{R}}(\mathit{\Lambda}^pT_x^\ast M,\mathit{\Lambda}^pT_y^\ast M)$, the homomorphisms between the vector spaces $\mathit{\Lambda}^pT_x^\ast M$ and $\mathit{\Lambda}^pT_y^\ast M$. We claim that the following two conditions determine the double forms $u^{i,p}(x,y)$ uniquely, therefore $H_N^p(t,x,y)$ are well-defined with these restrictions:
\begin{itemize}
\item $u^{0,p}(x,x)=\textrm{Id}\in\textrm{End}_{\mathbb{R}}(\mathit{\Lambda}^pT_x^\ast M)$, i.e., the identity endomorphsim of the vector space $\mathit{\Lambda}^pT_x^\ast M$;
\item \begin{equation}\label{eq:par}(\frac{\partial}{\partial t}+\Delta_y^p)H_N^p(t,x,y)=(4\pi t)^{-\frac{d}{2}}e^{-\frac{r^2(x,y)}{4t}}t^N\Delta_y^pu^{N,p}(x,y),\end{equation}
    where the superscript $p$ means the Laplacian acts on $p$-forms.
\end{itemize}
\textsc{Proof of the Claim}\\
Fix an $x\in M$ and introduce the Riemannian normal coordinates in the open neighborhood $U$ of $x$, then $g_{ij}(y)=\delta_{ij}+O(r^2(x,y))$. For $h(r(x,y))\in C^{\infty}(U)$ and $\omega\in C^{\infty}(U,\mathit{\Lambda}^pT^\ast M)$, computation gives
\begin{equation}\label{eq:del}\Delta_y^p(h\omega)=-(\frac{\partial^2h}{\partial r^2}+\frac{d-1}{r}\frac{\partial h}{\partial r}+\frac{1}{2g}\frac{\partial g}{\partial r}\frac{\partial h}{\partial r})\omega-\frac{2}{r}\frac{\partial h}{\partial r}\nabla_{r\frac{\partial}{\partial r}}\omega-h\Delta_y^p\omega,\end{equation}
where $g=\det(g_{ij})$ and $\nabla_{r\frac{\partial}{\partial r}}$ is the covariant derivative induced by the Levi-Civita connection with respect to $r\frac{\partial}{\partial r}$. Let $h=e^{-\frac{r^2(x,y)}{4t}}$ in (\ref{eq:del}) to get
\begin{equation}\Delta_y^p(e^{-\frac{r^2}{4t}}\omega)=-e^{-\frac{r^2}{4t}}\Big((\frac{r^2-2d}{4t}-\frac{r}{4gt}\frac{\partial g}{\partial r})\omega-\frac{1}{t}\nabla_{r\frac{\partial}{\partial r}}\omega-\Delta_y^p\omega\Big).\end{equation}
Therefore
\begin{eqnarray}
(\frac{\partial}{\partial t}+\Delta_y^p)H_N^p(t,x,y)=(4\pi t)^{-\frac{d}{2}}e^{-\frac{r^2}{4t}}\cdot\qquad\qquad\qquad\qquad\qquad\qquad\nonumber\\
\sum_{i=0}^{N}\Big((\frac{i}{t}+\frac{r}{4gt}\frac{\partial g}{\partial r})t^iu^{i,p}(x,y)+t^{i-1}\nabla_{r\frac{\partial}{\partial r}}u^{i,p}(x,y)+t^i\Delta_y^pu^{i,p}(x,y)\Big).
\end{eqnarray}
By (\ref{eq:par}), the coefficient of $(4\pi t)^{-\frac{d}{2}}e^{-\frac{r^2}{4t}}t^{i-1}$ in $(\frac{\partial}{\partial t}+\Delta_y^p)H_N^p(t,x,y)$ must be 0, it follows that
\begin{equation}\label{eq:pre}\nabla_{r\frac{\partial}{\partial r}}u^{i,p}(x,y)+(i+\frac{r}{4g}\frac{\partial g}{\partial r})u^{i,p}(x,y)=-\Delta_y^pu^{i-1,p}(x,y).\end{equation}
Now it is clear that to prove the claim, we only need to show that for every $\eta\in\mathit{\Lambda}^pT_x^\ast M$, the differential equations
\begin{equation}\label{eq:uip}\nabla_{r\frac{\partial}{\partial r}}u^{i,p}(\eta,y)+(i+\frac{r}{4g}\frac{\partial g}{\partial r})u^{i,p}(\eta,y)=-\Delta_y^pu^{i-1,p}(\eta,y),0\leq i\leq N\end{equation}
have unique solutions with the initial condition $u^{0,p}(\eta,x)=\eta$, here we adopt the convention $u^{-1,p}(x,y)\equiv0$.\\
Since we will argue by induction, it is convenient to rewrite equation (\ref{eq:uip}) in the following way:
\begin{equation}\label{eq:ri}\nabla_{r\frac{\partial}{\partial r}}\Big(r^ig^{\frac{1}{4}}u^{i,p}(\eta,y)\Big)=-r^ig^{\frac{1}{4}}\Delta_y^pu^{i-1,p}(\eta,y).\end{equation}
Fix a $y\in U$ and let $y(s),0\leq s\leq r(x,y)$ be the geodesic from $y$ to $x$. Use $\parallel_{y(s)}$ to denote the isomorphism $\mathit{\Lambda}^pT_y^\ast M\cong\mathit{\Lambda}^pT_{y(s)}^\ast M$ induced by the parallel translation along the geodesic. Let $u^{0,p}(\eta,y)=g^{-\frac{1}{4}}(y)$, then $u^{0,p}(\eta,x)=\eta$, equation (\ref{eq:ri}) is satisfied for $i=0$. Fix a $k\in\mathbb{N}_+$ and suppose that for $i<k$ we have determined the forms $u^{i,p}(\eta,y)$ satisfying (\ref{eq:ri}), then we define $u^{m,p}(\eta,y)$ as
\begin{eqnarray}
u^{m,p}(\eta,y)=-r(x,y)^{-m}g^{-\frac{1}{4}}(y)\cdot\qquad\qquad\qquad\qquad\qquad\qquad\qquad\nonumber\\
\int_0^{r(x,y)}\bigg(r\Big(x,y(s)\Big)\bigg)^{m-1}g^{\frac{1}{4}}\Big(y(s)\Big)\parallel_{y(s)}\bigg(\Delta_y^pu^{m-1,p}\Big(\eta,y(s)\Big)\bigg)\textrm{d}s.
\end{eqnarray}
It is clear that $u^{m,p}(\eta,y)\in C^{\infty}(U,\mathit{\Lambda}^pT^\ast M)$ and it can be checked that $u^{m,p}(\eta,y)$ satisfies the equation (\ref{eq:ri}) for $i=m$. This proves the existence part.\\
To prove the uniqueness, we derive from (\ref{eq:par}) the important relation
\begin{equation}iu^{i,p}(\eta,x)=-\Big(\Delta_y^pu^{i-1,p}(\eta,y)\Big)(\eta,x).\end{equation}
Compared with (\ref{eq:ri}), it suffices to prove that any $\zeta\in C^{\infty}(U,\mathit{\Lambda}^pT^\ast M)$ satisfying $\nabla_{r\frac{\partial}{\partial r}}\zeta=0$ and $\zeta(x)=0$ must be a zero section. This is trivially true because $\zeta(y)$ is just the parallel translation of $\zeta(x)$ along the geodesic for any point $y\in U$.\qquad$\square$\\
We can now construct the parametrix from $H_N^p(t,x,y)$. Let $W_1\times W_1$ and $W_2\times W_2$ be open sets such that $\overline{W_1\times W_1}\subset\overline{W_2\times W_2}\subset U\times U$. Choose $\phi(x,y)\in C^{\infty}(M\times M)$ such that $\phi(x,y)=1$ in $W_1\times W_1$ and $\phi(x,y)=0$ outside $W_2\times W_2$. We define
\begin{equation}G_N^p(t,x,y)=\phi(x,y)H_N^p(t,x,y),\end{equation}
and
\begin{equation}K_N^p(t,x,y)=(\frac{\partial}{\partial t}+\Delta_y^p)G_N^p(t,x,y).\end{equation}
\textbf{Lemma 5.2.4.} $G_N^p(t,x,y)$ \textit{is a parametrix for differential} $p$-\textit{forms when} $N>\frac{d}{2}$.\\
\textit{Proof}. We first show that $(\frac{\partial}{\partial t}+\Delta_y^p)G_N^p(t,x,y)$ is continuous at $t=0$. Since $G_N^p=0$ on $M\times M\setminus W_2\times W_2$, we only need to consider the problem on $\mathbb{R}_+\times W_2\times W_2$.\\
On $\mathbb{R}_+\times W_1\times W_1$, by (\ref{eq:par}) and the assumption $N>\frac{d}{2}$
\begin{eqnarray}
\lim_{t\rightarrow 0}(\frac{\partial}{\partial t}+\Delta_y^p)G_N^p(t,x,y)=\lim_{t\rightarrow 0}(\frac{\partial}{\partial t}+\Delta_y^p)H_N^p(t,x,y)\qquad\qquad\nonumber\\
=\lim_{t\rightarrow 0}(4\pi t)^{-\frac{d}{2}}e^{-\frac{r^2(x,y)}{4t}}t^N\Delta_y^pu^{N,p}(x,y)\nonumber\\
=0.\qquad\qquad\qquad\qquad\qquad\qquad\qquad
\end{eqnarray}
Finally, on $W_2\times W_2\setminus W_1\times W_1$, use (\ref{eq:del}) we have
\begin{eqnarray}
\label{eq:so}(\frac{\partial}{\partial t}+\Delta_y^p)G_N^p(t,x,y)=\lim_{t\rightarrow 0}(\frac{\partial}{\partial t}+\Delta_y^p)\Big(\phi H_N^p(t,x,y)\Big)\qquad\qquad\qquad\nonumber\\
=-(\frac{\partial^2\phi}{\partial r^2}+\frac{d-1}{r}\frac{\partial \phi}{\partial r}+\frac{1}{2g}\frac{\partial g}{\partial r}\frac{\partial \phi}{\partial r})H_N^p(t,x,y)\nonumber\\
-\frac{2}{r}\frac{\partial \phi}{\partial r}\nabla_{r\frac{\partial}{\partial r}}H_N^p(t,x,y)-\phi\Delta_y^pH_N^p(t,x,y)\quad\nonumber\\
=\psi_1(x,y)H_N^p+\psi_2(x,y)\nabla_{r\frac{\partial}{\partial r}}H_N^p-\phi\Delta_y^pH_N^p,
\end{eqnarray}
where $\psi_1,\psi_2\in C^\infty(M\times M)$. We can thus deduce from (\ref{eq:so}) that
\begin{equation}\lim_{t\rightarrow 0}(\frac{\partial}{\partial t}+\Delta_y^p)G_N^p(t,x,y)=0.\end{equation}
To show that $G_N^p(t,x,y)$ satisfies the second condition in Definition 5.2.3, we may assume that there exists an $\varepsilon>0$ such that $B_\varepsilon(x)\subset W_1$, where $B_\varepsilon(x)$ is the open ball centering at $x$ with radius $\varepsilon$. We then have
\begin{eqnarray}
\lim_{t\rightarrow 0}\int_M(4\pi t)^{-\frac{d}{2}}\phi(x,y)e^{-\frac{r^2}{4t}}u^{i,p}(x,y)\wedge\star\omega(y)\quad\qquad\nonumber\\
=\lim_{t\rightarrow 0}\int_{B_\varepsilon(x)}(4\pi t)^{-\frac{d}{2}}\phi(x,y)e^{-\frac{r^2}{4t}}u^{i,p}(x,y)\wedge\star\omega(y)\qquad\nonumber\\
=\lim_{t\rightarrow 0}\int_{M\setminus B_\varepsilon(x)}(4\pi t)^{-\frac{d}{2}}\phi(x,y)e^{-\frac{r^2}{4t}}u^{i,p}(x,y)\wedge\star\omega(y).\quad
\end{eqnarray}
Since $r\geq\varepsilon$, the second integral above vanishes. To compute the first integral, we use the exponential map $\exp_x$ to pull the integral back to $T_xM$:
\begin{eqnarray}
\int_{B_\varepsilon(x)}(4\pi t)^{-\frac{d}{2}}\phi(x,y)e^{-\frac{r^2}{4t}}u^{i,p}(x,y)\wedge\star\omega(y)\qquad\qquad\qquad\qquad\qquad\qquad\quad\nonumber\\
=\int_{B_\varepsilon(0)\subset T_xM}J(\exp_x)(4\pi t)^{-\frac{d}{2}}\phi(x,y)e^{-\frac{r^2(0,v)}{4t}}u^{i,p}(x,\exp_xv)\wedge\star\omega(\exp_xv)\quad\nonumber\\
=\int_{T_xM}J(\exp_x)(4\pi t)^{-\frac{d}{2}}\phi(x,y)e^{-\frac{r^2(0,v)}{4t}}u^{i,p}(x,\exp_xv)\wedge\star\omega(\exp_xv),\qquad\quad
\end{eqnarray}
where $\exp_x(v)=y$, $u^{i,p}$ are extended to be 0 outside $B_\varepsilon(x)$ and $J(\exp_x)$ denotes the Jacobian of $\exp_x$.\\
Since $(4\pi t)^{-\frac{d}{2}}e^{-\frac{r^2}{4t}}$ is the ordinary heat kernel of $\mathbb{R}^d$, by Definition 5.2.2 we have
\begin{eqnarray}
\lim_{t\rightarrow 0}\int_{T_xM\cong\mathbb{R}^d}J(\exp_x)(4\pi t)^{-\frac{d}{2}}\phi(x,y)e^{-\frac{r^2(0,v)}{4t}}u^{i,p}(x,\exp_xv)\wedge\star\omega(\exp_xv)\nonumber\\
=\star\Big(u^{i,p}(x,x)\wedge\star\omega(x)\Big)=\Big\langle u^{i,p}(x,x),\omega(x)\Big\rangle.\qquad\qquad\qquad\qquad\qquad\quad
\end{eqnarray}
Since $u^{0,p}(x,x)=\textrm{Id}$ by our construction, we finally get
\begin{equation}\label{eq:fige}\lim_{t\rightarrow 0}\int_M(4\pi t)^{-\frac{d}{2}}e^{-\frac{r^2}{4t}}\phi(x,y)\Big(\sum_{i=0}^{N}t^iu^{i,p}(x,y)\Big)\wedge\star\omega(y)=\omega(x),\end{equation}
since all the summands in (\ref{eq:fige}) vanishes trivially except the first one.\qquad$\square$\\
We now complete the construction of the heat kernel. For a metric $\langle\cdot,\cdot\rangle_z$ on $\mathit{\Lambda}^pT^\ast M$, there is a natural way to extend it to $\mathit{\Lambda}^pT^\ast M\otimes\mathit{\Lambda}^pT^\ast M$:
\begin{equation}\Big\langle\eta(x)\otimes\omega_1(z),\zeta(y)\otimes\omega_2(z)\Big\rangle=\Big\langle\omega_1(z),\omega_2(z)\Big\rangle_z\eta(x)\otimes\zeta(y),\end{equation}
where $\eta(x)\otimes\omega_1(z),\zeta(y)\otimes\omega_2(z)\in C^{\infty}(M\times M,\mathit{\Lambda}^pT^\ast M\otimes\mathit{\Lambda}^p T^\ast M)$. Using this notation, we define
\begin{eqnarray}K^0(t,x,y)=K_N^p(t,x,y),\qquad\qquad\qquad\qquad\qquad\qquad\qquad\qquad\qquad\nonumber\\
K^m(t,x,y)=\int_0^t\textrm{d}s\int_M\Big\langle K^{m-1}(s,x,z),K_N^p(t-s,z,y)\Big\rangle\textrm{dvol}_z,m\geq 1,
\end{eqnarray}
and
\begin{eqnarray}\label{eq:ep}e^p(t,x,y)=G_N^p(t,x,y)+\qquad\qquad\qquad\qquad\qquad\qquad\qquad\qquad\qquad\nonumber\\
\sum_{m\geq 0}(-1)^{m+1}\int_0^t\textrm{d}s\int_M\Big\langle K^m(s,x,z),G_N^p(t-s,z,y)\Big\rangle\textrm{dvol}_z.
\end{eqnarray}
\textbf{Theorem 5.2.5.} $e^p(t,x,y)$ \textit{defined above is the heat kernel for differential} $p$-\textit{forms}.\\
\textit{Proof}. We first show that $e^p(t,x,y)$ is well defined, i.e., the right hand side of (\ref{eq:ep}) converges everywhere to a smooth double form.\\
We do some localization to make the coordinates available in our estimate. Since $M$ is compact, there is a finite family of open sets $\{W_i\}_{i=1}^{n}$ such that $M=\bigcup_{i=1}^{n}W_i$ and $W_i\subset U_i$, where $\{U_i\}_{i=1}^{n}$ is the coordinate chart of $M$. Let $\{\varphi_i\}_{i=1}^n$ be a partition of unity relative to the covering $\{W_i\}_{i=1}^{n}$, and let $\{\lambda_i\}_{i=1}^n$ be $C^\infty$-functions with  supports in $U_i$ which are identically 1 on $W_i$.\\
For any $\xi(x,y)\in C^\infty(M\times M,\mathit{\Lambda}^pT^\ast M\otimes\mathit{\Lambda}^pT^\ast M)$, we define the localized form $\xi_{i,j}(x,y)\in C^\infty(U_i\times U_j,\mathit{\Lambda}^pT^\ast M\otimes\mathit{\Lambda}^pT^\ast M)$ as
\begin{equation}\xi_{i,j}(x,y)=\lambda_i(x)\lambda_j(y)\xi(x,y).\end{equation}
As a consequence, any $\mu(x,y)\in C^\infty(U_i\times U_j,\mathit{\Lambda}^pT^\ast M\otimes\mathit{\Lambda}^pT^\ast M)$ can be expressed by
\begin{equation}\mu(x,y)=a_{IJ}(x,y)dx^{I}\otimes dy^{J},\end{equation}
where $I=\{i_k|1\leq k\leq p,i_1<\cdot\cdot\cdot<i_p\}$ and $J=\{j_k|1\leq k\leq p,j_1<\cdot\cdot\cdot<j_p\}$ are multi-indices and $a_{IJ}(x,y)\in C^\infty(U_i\times U_j)$. To begin our estimate, it is natural to introduce the following norm:
\begin{equation}\Big\|\mu(x,y)\Big\|_{i,j}=\sum_{I,J}\sup_{x\in U_i,y\in U_j}\Big|a_{IJ}(x,y)\Big|,\end{equation}
It is clear that to bound the right hand side of (\ref{eq:ep}), we only need to bound $K^m$, i.e., to bound the norm $\|K^m_{i,j}\|_{i,j}$. By (\ref{eq:par}), there exists a constant $C_1$ such that $\|K_{i,j}^0\|_{i,j}\leq C_1t^{N-\frac{d}{2}}$. We can then proceed by induction to obtain a bound for $\|K^m_{i,j}\|_{i,j}$. Suppose that we have
\begin{equation}\Big\|K^{m-1}_{i,j}(t,x,y)\Big\|_{i,j}\leq (C_2C_1)^mt^{m(N-\frac{d}{2})+m-1}\frac{\Gamma^m(N-\frac{d}{2}+1)}{\Gamma\Big(m(N-\frac{d}{2})+m\Big)}\end{equation}
for some constant $C_2$. By the definition of $K_{i,j}^m$, we have
\begin{eqnarray}K^m_{i,j}(t,x,y)=\qquad\qquad\qquad\qquad\qquad\qquad\qquad\qquad\qquad\qquad\qquad\nonumber\\
\int_0^t\textrm{d}s\int_M\Big\langle\lambda_i(x)\sum_{i=1}^{n}\varphi_i(z) K^{m-1}(s,x,z),\lambda_j(y)K_N^p(t-s,z,y)\Big\rangle\textrm{dvol}_z.
\end{eqnarray}
It follows easily that
\begin{eqnarray}
\Big\|K^m_{i,j}(t,x,y)\Big\|_{i,j}\leq C_3C_1(C_2C_1)^m\frac{\Gamma^m(N-\frac{d}{2}+1)}{\Gamma\Big(m(N-\frac{d}{2})+m\Big)}\cdot\nonumber\\
\int_0^ts^{m(N-\frac{d}{2})+m-1}(t-s)^{N-\frac{d}{2}}\textrm{d}s,\qquad
\end{eqnarray}
where $C_3$ is a constant independent of $m$. Since we may assume that $C_2>C_3$ and $C_1\geq 1$, we obtain
\begin{equation}\label{eq:km}\Big\|K^{m}_{i,j}(t,x,y)\Big\|_{i,j}\leq (C_2C_1)^{m+1}t^{(m+1)(N-\frac{d}{2})+m}\frac{\Gamma^{m+1}(N-\frac{d}{2}+1)}{\Gamma\Big((m+1)(N-\frac{d}{2})+m+1\Big)}.\end{equation}
This shows that the right hand side of (\ref{eq:ep}) converges uniformly to a smooth double $p$-form. Since
\begin{eqnarray}
(\frac{\partial}{\partial t}+\Delta_y^p)e^p(t,x,y)\qquad\qquad\qquad\qquad\qquad\qquad\qquad\qquad\nonumber\\
=\sum_{m=0}^\infty(-1)^{m+1}\Big(K^m(t,x,y)+K^{m+1}(t,x,y)\Big)+K_N^p(t,x,y)\nonumber\\
=K_N^p(t,x,y)-K_N^p(t,x,y)=0,\qquad\qquad\qquad\qquad\qquad\quad
\end{eqnarray}
and
\begin{eqnarray}
\lim_{t\rightarrow 0}\int_Me^p(t,x,y)\wedge\star\omega(y)\qquad\qquad\qquad\qquad\qquad\qquad\qquad\qquad\qquad\qquad\quad\nonumber\\
=\lim_{t\rightarrow 0}\Bigg(\int_MG_N^p(t,x,y)\wedge\star\omega(y)+\qquad\qquad\qquad\qquad\qquad\qquad\qquad\qquad\qquad\nonumber\\
\sum_{m=0}^\infty(-1)^{m+1}\int_0^t\textrm{d}s\int_M\bigg(\int_M\Big\langle K^m(s,x,z),G_N^p(t-s,z,y)\Big\rangle\textrm{dvol}_z\bigg)\wedge\star\omega(y)\Bigg)\nonumber\\
=\omega(x)\quad\qquad\qquad\qquad\qquad\qquad\qquad\qquad\qquad\qquad\qquad\qquad\qquad\qquad\qquad\nonumber
\end{eqnarray}
by Lemma 5.2.4 and the estimate (\ref{eq:km}), $e^p(t,x,y)$ is the heat kernel by Definition 5.2.2.\qquad$\square$\\
We finish this subsection by expressing $e^p(t,x,y)$ with respect to a basis of $L^2(M,\mathit{\Lambda}^pT^\ast M\otimes\mathit{\Lambda}^pT^\ast M)$. This expression will lead to an integral expression of the trace of the heat operator, which will be used in our discussion in the next section.\\
\textbf{Proposition 5.2.6.} \textit{Let} $\{\beta^{i,p}\}_{i\in\mathbb{Z}_+}$ \textit{be an orthonormal basis of} $L^2(M,\mathit{\Lambda}^pT^\ast M)$ \textit{satisfying} $\Delta^p\beta^{i,p}=\lambda_i^p\beta^{i,p}$. \textit{Then we have the pointwise convergence}
\begin{equation}\label{eq:ee}e^p(t,x,y)=\sum_{i=0}^\infty e^{-\lambda_i^pt}\beta^{i,p}(x)\otimes\beta^{i,p}(y).\end{equation}
\textit{Proof}. Since $e^p(t,x,y)\in L^2(\mathbb{R}_+\times M\times M,\mathit{\Lambda}^pT^\ast M\otimes\mathit{\Lambda}^pT^\ast M)$, we have
\begin{equation}\label{eq:ebeta}e^p(t,x,y)=e_i^p(t,x)\beta^{i,p}(y),e_i^p(t,x)=\star\Big(e^p(t,x,y)\wedge\star\beta^{i,p}(y)\Big).\end{equation}
Then by the self-adjoint property of $\Delta_y^p$:
\begin{eqnarray}
\frac{\partial}{\partial t}e_i^p(t,x)=\frac{\partial}{\partial t}\star\Big(e^p(t,x,y)\wedge\star\beta^{i,p}(y)\Big)\quad\nonumber\\
=-\star\Big(\Delta_y^pe^p(t,x,y)\wedge\star\beta^{i,p}(y)\Big)\nonumber\\
=-\star\Big(e^p(t,x,y)\wedge\star\Delta_y^p\beta^{i,p}(y)\Big)\nonumber\\
=-\lambda_i^p\star\Big(e^p(t,x,y)\wedge\star\beta^{i,p}(y)\Big)\nonumber\\
\label{eq:beta}=-\lambda_i^pe^{i,p}(t,x).\qquad\qquad\qquad\quad
\end{eqnarray}
Solving the differential equation (\ref{eq:beta}) to get $e_i^p(t,x)=c_i^p(x)e^{-\lambda_i^p t}$.\\
For an arbitrary $\omega=a_i\beta^{i,p}\in L^2(\mathit{\Lambda}^pT^\ast M)$ with $a_i\in\mathbb{R}$, we have by Definition 5.2.2 that
\begin{eqnarray}
\omega(x)=\lim_{t\rightarrow 0}\star\Big(e^p(t,x,y)\wedge\star\omega(y)\Big)\qquad\qquad\qquad\quad\nonumber\\
=\lim_{t\rightarrow 0}\star\Big(\sum_{i=0}^\infty e^{-\lambda_i^pt}c_i^p(x)\beta^{i,p}(y)\wedge\star a_j\beta^{j,p}(y)\Big)\nonumber\\
=\lim_{t\rightarrow 0}\sum_{i=0}^\infty e^{-\lambda_i^pt}c_i^p(x)a_i=\sum_{i=0}^\infty c_i^p(x)a_i,\qquad\quad
\end{eqnarray}
which implies $c_i^p(x)=\beta^{i,p}(x)$. Then (\ref{eq:ebeta}) becomes
\begin{equation}e^p(t,x,y)=\sum_{i=0}^\infty e^{-\lambda_i^p t}\beta^{i,p}(x)\otimes\beta^{i,p}(y)\in L^2(\mathbb{R}_+\times M\times M,\mathit{\Lambda}^pT^\ast M\otimes\mathit{\Lambda}^pT^\ast M).\end{equation}
As a result, there exists a subsequence of the right hand side of (\ref{eq:ee}) which converges a.e. to $e^p(t,x,y)$ on $M$.\\
By Parseval's equality
\begin{equation}\Big\langle e^p(\frac{t}{2},x,z),e^p(\frac{t}{2},y,z)\Big\rangle_z=\sum_{i=0}^\infty e^{-\lambda_i^pt}\beta^{i,p}(x)\otimes\beta^{i,p}(y).\end{equation}
This shows that the right hand side of (\ref{eq:ee}) converges everywhere to a continuous limit, and the limit must be $e^p(t,x,y)$ by the discussions above.\qquad$\square$\\
\textbf{Corollary 5.2.7.}
\begin{equation}\sum_{i=0}^\infty e^{-\lambda_i^pt}=\textrm{Tr}(\frac{\partial}{\partial t}+\Delta^p)=\int_M\textrm{tr}_xe^p(t,x,x)\textrm{dvol}_x,\end{equation}
\textit{where} $\textrm{Tr}$ \textit{denotes the trace of an operator and} $\textrm{tr}_x$ \textit{denotes the trace of} $e^p(t,y,x)$ \textit{at} $y=x$.\\
\textit{Proof}. By the above proposition,
\begin{eqnarray}
\int_M\textrm{tr}_xe^p(t,x,x)\textrm{dvol}_x=\int_M\textrm{tr}_x\Big(\sum_{i=0}^\infty e^{-\lambda_i^pt}\beta^{i,p}(x)\otimes\beta^{i,p}(x)\Big)\textrm{dvol}_x\nonumber\\
=\sum_{i=0}^\infty e^{-\lambda_i^pt}\int_M\textrm{tr}_x\Big(\beta^{i,p}(x)\otimes\beta^{i,p}(x)\Big)\textrm{dvol}_x\nonumber\\
=\sum_{i=0}^\infty e^{-\lambda_i^pt}\Big\langle\beta^{i,p}(x),\beta^{i,p}(x)\Big\rangle_x\quad\qquad\qquad\nonumber\\
=\sum_{i=0}^\infty e^{-\lambda_i^pt}=\textrm{Tr}(\frac{\partial}{\partial t}+\Delta_y^p).\qquad\qquad\qquad
\end{eqnarray}
The proof is complete.\qquad$\square$\\
\textsc{Remark}\\
This corollary justifies the abbreviation $e^{-t\Delta^p}=\frac{\partial}{\partial t}+\Delta^p$. We will assume this convention in the following discussions.
\subsection{Supertrace and the McKean-Singer Conjecture}
To describe the McKean-Singer conjecture we need the asymptotic expansion of the trace of $e^{-t\Delta^p}$.\\
\textbf{Definition 5.3.1.} For $t\in\mathbb{R}$, we say that the formal power series $\sum_{i=i_0}^\infty a_it^i$ is the \textit{asymptotic expansion} of $A(t)\in C^\infty(\mathbb{R})$ and write $A(t)\sim\sum_{i=i_0}^\infty a_it^i$ if for all $N\geq i_0$, we have
\begin{equation}\lim_{t\rightarrow 0}\frac{A(t)-\sum_{i=i_0}^Na_it^i}{t^N}=0.\end{equation}
\textbf{Proposition 5.3.2.} \textit{Suppose} $\{\lambda_i^p\}$ \textit{is the spectrum of} $\Delta_x^p$, \textit{then}
\begin{equation}\sum_{i=0}^\infty e^{\lambda_i^pt}\sim(4\pi t)^{\frac{d}{2}}\sum_{i=0}^\infty\Big(\int_M\textrm{tr}_xu^{i,p}(x,x)\textrm{dvol}_x\Big)t^i.\end{equation}
\textit{Proof}. It is straightforward from (\ref{eq:hnp}) and (\ref{eq:ep}) that
\begin{equation}e^p(t,x,x)\sim(4\pi t)^{-\frac{d}{2}}\sum_{i=0}^{\infty}u^{i,p}(x,x)t^i.\end{equation}
By Corollary 5.2.7, we have
\begin{eqnarray}
\sum_{i=0}^\infty e^{-\lambda_i^pt}=\int_M\textrm{tr}_xe^p(t,x,x)\textrm{dvol}_x\sim(4\pi t)^{-\frac{d}{2}}\int_M\sum_{i=0}^\infty u^{i,p}(x,x)t^i\textrm{dvol}_x\nonumber\\
\sim(4\pi t)^{\frac{d}{2}}\sum_{i=0}^\infty\Big(\int_M\textrm{tr}_xu^{i,p}(x,x)\textrm{dvol}_x\Big)t^i.\qquad\qquad\qquad\qquad
\end{eqnarray}
This finishes the proof.\qquad$\square$\\
We begin with the following crucial observation of McKean and Singer (\cite{ms}):\\
\textbf{Lemma 5.3.3.} \textit{Let} $\Delta^p$ \textit{be the Laplacian on} $p$-\textit{forms}. \textit{For} $\lambda\in\mathbb{R}_+$, \textit{let} $E^p_\lambda$ \textit{be the} $\lambda$-\textit{eigenspace of} $\Delta^p$. \textit{Then we have the following exact sequence}:
\begin{equation}0\rightarrow E_\lambda^0\xrightarrow{d}\cdot\cdot\cdot\xrightarrow{d}E_\lambda^d\rightarrow 0.\end{equation}
\textit{Proof}. First note that the sequence is well defined, since if $\omega\in E_\lambda^p$, then $\Delta^{p+1}d\omega=d\Delta^p\omega=\lambda d\omega$.\\
If $\omega\in E_\lambda^p$ has $d\omega=0$, then
\begin{equation}d(\frac{1}{\lambda}\delta\omega)=\frac{1}{\lambda}\Delta^p\omega-\frac{1}{\lambda}\delta d\omega=\frac{1}{\lambda}\Delta^p\omega=\omega,\end{equation}
this finishes the proof.\qquad$\square$\\
\textbf{Corollary 5.3.4.}
\begin{equation}\label{eq:dim}\sum_{p=0}^d(-1)^p\dim_{\mathbb{R}}E_\lambda^p=0,\lambda>0.\end{equation}
\textit{Proof}. By the above lemma, $d+\delta$ defines an isomorphism $\bigoplus_{k=0}^{\frac{d}{2}}E^{2k}_{\lambda}\rightarrow\bigoplus_{k=0}^{\frac{d}{2}-1}E^{2k+1}_{\lambda}$, then (\ref{eq:dim}) follows directly.\qquad$\square$\\
Recall that the spectrum of the Laplacian on a compact manifold must be discrete, then we have the following important corollary:\\
\textbf{Corollary 5.3.5.} \textit{Let} $\{\lambda_i\}$ \textit{be the spectrum of} $\Delta^p$. \textit{Then}
\begin{equation}\sum_{p=0}^{d}(-1)^p\sum_{i=0}^\infty e^{-\mathit{\lambda}_i^pt}=\sum_{p=0}^d\dim_{\mathbb{R}}\ker\Delta^p.\end{equation}
\textit{Proof}. By the above corollary, we have
\begin{equation}\sum_{p=0}^{d}(-1)^p\sum_{i=0}^\infty e^{-\lambda_i^pt}=\sum_{p=0}^d(-1)^p\sum_{i,\lambda_i^p=0}e^{-\lambda_i^pt}=\sum_{p=0}^d\dim_{\mathbb{R}}\ker\Delta^p.\end{equation}
This proves the corollary.\qquad$\square$\\
This leads naturally to the following definition.\\
\textbf{Definition 5.3.6.} Let $V$ be a vector bundle over $M$ and $V^\ast$ its dual bundle. For every $A\in\textrm{End}_{\mathbb{R}}(\mathit{\Lambda}^\ast V^\ast)$, we define the \textit{supertrace} $\textrm{Tr}\Big((-1)^FA\Big)$ of $A$ as the trace of $A$ on even forms minus the trace of $A$ on odd forms.\\
Using this definition, Corollary 5.3.5 reads $\textrm{Tr}\Big((-1)^Fe^{-t\Delta}\Big)=\sum_{p=0}^d\dim_{\mathbb{R}}\ker\Delta^p$. It leads directly to the following surprising result:\\
\textbf{Theorem 5.3.7.}
\begin{equation}\label{eq:msc}(4\pi)^{-\frac{d}{2}}\int_M\sum_{p=0}^d(-1)^p\textrm{tr}_xu^{i,p}(x,x)\textrm{dvol}_x=\left\{\begin{array}{ll}0 & i\neq\frac{d}{2}\\\chi(M) & i=\frac{d}{2},\dim M \textrm{ even}\end{array}\right.\end{equation}
\textit{Proof}. First note that $\textrm{Tr}\Big((-1)^Fe^{-t\Delta}\Big)$ is independent of $t$ by Corollary 5.3.5. By the Hodge theorem, we have
\begin{eqnarray}
\chi(M)=\sum_{p=0}^d(-1)^p\dim_{\mathbb{R}}H_{\textrm{dR}}^p(M,\mathbb{R})=\sum_{p=0}^d\dim_{\mathbb{R}}\ker\Delta^p\qquad\quad\nonumber\\
=\textrm{Tr}\Big((-1)^Fe^{-t\Delta}\Big)=\sum_{p=0}^d(-1)^p\int_M\textrm{tr}_xe^p(t,x,x)\textrm{dvol}_x,\quad
\end{eqnarray}
Using Proposition 5.3.2, we get the asymptotic expansion of $\chi(M)$:
\begin{equation}\label{eq:chi}\chi(M)\sim(4\pi t)^{-\frac{d}{2}}\sum_{i=0}^\infty\Big(\int_M\sum_{p=0}^d(-1)^p\textrm{tr}_xu^{i,p}(x,x)\textrm{dvol}_x\Big)t^i.\end{equation}
Since $\chi(M)$ is independent of $t$, only the constant term on the right hand side of (\ref{eq:chi}) can be nonzero.\qquad$\square$\\
Compare this theorem with the Gauss-Bonnet-Chern theorem (\ref{eq:gbco}), it is reasonable to make the following conjecture.\\
\textbf{Conjecture 5.3.8. (McKean-Singer)} \textit{For every even} $d$,
\begin{equation}\Omega=(4\pi)^{-\frac{d}{2}}\sum_{p=1}^d(-1)^p\textrm{tr}_xu^{\frac{d}{2},p}(x,x)\textrm{dvol}_x,\end{equation}
\textit{where} $\Omega$ \textit{is the Gauss-Bonnet integrand defined by} (\ref{eq:gbi}).\\
\textsc{Remarks}
\begin{itemize}
\item The conjecture is trivial for $d=2$, since it is well-known that $u^{1,0}(x,x)=\frac{1}{6}R_i^i$, where $R_i^i$ is the scalar curvature. This leads to a proof of the classical Gauss-Bonnet theorem on surfaces. In \cite{ms}, McKean and Singer also proved their conjecture for $d=4$ by expressing $u^{2,0}(x,x)$ in terms of curvature tensors. Although it is known that for any $i>0$, $u^{i,0}(x,x)$ is a sum of universal polynomials of curvature tensors, the computation for $i\geq 3$ seems hopelessly complicated. The interested reader may refer to \cite{bg}, \cite{pg} and \cite{sr} for these results.
\item Note that by combining the local property (asymptotic expansion of the heat kernel) with the global property (spectral theory of the Laplacian), we have successfully related curvature with the topological invariant $\chi(M)$.
\end{itemize}
\subsection{Proof of the McKean-Singer Conjecture}
It is obvious that to complete our proof of the Gauss-Bonnet-Chern theorem, one only need to prove the McKean-Singer conjecture. This is first done by Patodi in \cite{vp}, using classical tensor calculus. We will present his proof here since this is the most direct and self-contained approach. Another proof using invariance theory can be found in \cite{pg}. A third proof using fermion calculus is contained in \cite{sr}. A fourth proof using calculus of Clifford algebras was discovered in \cite{eg}.\\
Since our proof is based on tensor calculus, we restate our main result in this section as follows:\\
\textbf{Theorem 5.4.1. (McKean-Singer-Patodi)} \textit{The following formula holds}:
\begin{eqnarray}(4\pi)^{-\frac{d}{2}}\sum_{p=0}^d(-1)^p\textrm{tr}_xu^{\frac{d}{2},p}(x,x)=\quad\qquad\qquad\qquad\qquad\qquad\qquad\qquad\nonumber\\
c_d\sum_{\sigma_1,\sigma_2\in\Sigma_{d}}\textrm{sgn }\sigma_1\textrm{sgn }\sigma_2R_{\sigma_1(1)\sigma_1(2)\sigma_2(1)\sigma_2(2)}\cdot\cdot\cdot R_{\sigma_1(d-1)\sigma_1(d)\sigma_2(d-1)\sigma_2(d)},
\end{eqnarray}
\textit{where} $c_d=\frac{(-1)^{\frac{d}{2}}}{(8\pi)^{\frac{d}{2}}(\frac{d}{2})!}$.\\
As we shall see, the proof of the McKean-Singer conjecture is a process of canceling terms. The remarkable cancelation lemma proved by Patodi in \cite{vp} actually provides much more than Theorem 5.4.1 requires. To prove Patodi's cancelation lemma, a series of technical lemmas will be needed, and we will establish them in the sequel. We first introduce two useful operators which will appear frequently in our computations.\\
Let $V$ be a vector bundle over $M$ of rank $d$, $V^\ast$ its dual. Then we can form the exterior product $\mathit{\Lambda}^pV$. An operator $\mathcal{A}\in\textrm{End}_{\mathbb{R}}(V)$ extends naturally to an element of $\mathit{\Lambda}^pV$ in the following two ways:
\begin{equation}\Lambda^p\mathcal{A}(v_1\wedge\cdot\cdot\cdot\wedge v_p)=\mathcal{A}v_1\wedge\cdot\cdot\cdot\wedge \mathcal{A}v_p,v_i\in V,1\leq i\leq p;\end{equation}
\begin{equation}\textrm{D}^p\mathcal{A}(v_1\wedge\cdot\cdot\cdot\wedge v_p)=\sum_{i=1}^pv_1\wedge\cdot\cdot\cdot\wedge v_{i-1}\wedge\mathcal{A}v_i\wedge v_{i+1}\wedge\cdot\cdot\cdot\wedge v_p,\end{equation}
where we have adopted the convention that $\Lambda^0\mathcal{A}=\textrm{Id}\in\textrm{End}(\mathbb{R})$ and $\textrm{D}^0\mathcal{A}=0\in\textrm{End}(\mathbb{R})$.\\
The following three lemmas are purely algebraic and do not include covariant derivatives, they
can be regarded as a baby version of Patodi's cancelation lemma and its applications.\\
\textbf{Lemma 5.4.2.} \textit{Let} $\mathcal{A}_1,\cdot\cdot\cdot,\mathcal{A}_k\in\textrm{End}_{\mathbb{R}}(V)$, $k\leq d$. \textit{When} $k=d$, \textit{suppose that} $\det(x_1\mathcal{A}_1+\cdot\cdot\cdot+x_d\mathcal{A}_d)=a_1x_1^d+\cdot\cdot\cdot+a_dx_d^d+\cdot\cdot\cdot+a_{12\cdot\cdot\cdot d}x_1\cdot\cdot\cdot x_d$ \textit{for some coefficients} $a_1,\cdot\cdot\cdot,a_d,\cdot\cdot\cdot,a_{12\cdot\cdot\cdot d}\in\mathbb{R}$. \textit{Then we have}
\begin{equation}\textrm{Tr}\Big((-1)^F\textrm{D}^p\mathcal{A}_1\circ\cdot\cdot\cdot\circ\textrm{D}^p\mathcal{A}_k\Big)=\left\{\begin{array}{ll}0 & k<d\\(-1)^da_{12\cdot\cdot\cdot d} & k=d\end{array}\right.\end{equation}
\textit{the notation} $\textrm{Tr}\Big((-1)^F\textrm{D}^p\mathcal{A}_1\circ\cdot\cdot\cdot\circ\textrm{D}^p\mathcal{A}_k\Big)$ \textit{is justified by the identification between} $V$ \textit{and} $V^\ast$.\\
\textit{Proof}. Elementary linear algebra yields
\begin{eqnarray}
\det(I-e^{x_1\mathcal{A}_1}\cdot\cdot\cdot e^{x_k\mathcal{A}_k})=\textrm{Tr}\Big((-1)^F\Lambda^p(e^{x_1\mathcal{A}_1}\cdot\cdot\cdot e^{x_k\mathcal{A}_k})\Big)\nonumber\\
=\textrm{Tr}\Big((-1)^Fe^{x_1\textrm{D}^p\mathcal{A}_1}\cdot\cdot\cdot e^{x_k\textrm{D}^p\mathcal{A}_k}\Big).
\end{eqnarray}
The proof will be completed by comparing the coefficients of $x_1\cdot\cdot\cdot x_k$ in the above equality.\qquad$\square$\\
The next lemma needs some explanation of the notations. Recall that the tensor product $V^\ast\otimes V$ can be identified with $\textrm{End}_{\mathbb{R}}(V)$, therefore we can introduce the map $\rho_p:C^\infty\Big(U,(V^\ast\otimes V)\times(V^\ast\otimes V)\Big)\rightarrow\textrm{End}_{\mathbb{R}}(\mathit{\Lambda}^pV)$ by $\rho_p(\mathcal{A},\mathcal{B})=\textrm{D}^p\mathcal{A}\circ\textrm{D}^p\mathcal{B}$. This is a bilinear map and will therefore induce a linear map $\widetilde{\rho}_p:C^\infty(U,V^\ast\otimes V\otimes V^\ast\otimes V)\rightarrow\textrm{End}_{\mathbb{R}}(\mathit{\Lambda}^pV)$. For $\mathcal{A}\in C^\infty(U,V^\ast\otimes V\otimes V^\ast\otimes V)$, we shall identify $\textrm{D}^p\mathcal{A}$ with $\widetilde{\rho}_p(\textrm{D}^p\mathcal{A})$.\\
\textbf{Lemma 5.4.3.} \textit{Let} $l\in\mathbb{Z}_+$ \textit{such that} $l<\frac{d}{2}$, $\sigma\in\Sigma_d$. $\mathcal{A}_1,\cdot\cdot\cdot,\mathcal{A}_l\in C^\infty(U,V^\ast\otimes V\otimes V^\ast\otimes V)$. \textit{Then}
\begin{equation}\textrm{Tr}\Big((-1)^F\textrm{D}^p\mathcal{A}_{\sigma(1)}\circ\cdot\cdot\cdot\circ\textrm{D}^p\mathcal{A}_{\sigma(l)}\Big)=0.\end{equation}
\textit{Proof}. This is an easy corollary of Lemma 5.4.2.\qquad$\square$\\
Let $\{e_1,\cdot\cdot\cdot,e_d\}$ be smooth sections of $V$ which form a local frame on $U$, and $\{e_1^\ast,\cdot\cdot\cdot,e_d^\ast\}$ its fiberwise dual, then $\mathcal{A}\in C^\infty(U,V^\ast\otimes V\otimes V^\ast\otimes V)$ has the form $\mathcal{A}=a^{ijkl}e_i^\ast\otimes e_j\otimes e_k^\ast\otimes e_l$, where $a^{ijkl}\in\mathbb{R}$.\\
\textbf{Lemma 5.4.4.} \textit{With the abbreviation} $(\textrm{D}^p\mathcal{A})^{\frac{d}{2}}=\textrm{D}^p\mathcal{A}\circ\cdot\cdot\cdot\circ\textrm{D}^p\mathcal{A}$, \textit{we have}
\begin{eqnarray}
\textrm{Tr}\Big((-1)^F(\textrm{D}^p\mathcal{A})^{\frac{d}{2}}\Big)=\quad\qquad\qquad\qquad\qquad\qquad\qquad\qquad\qquad\nonumber\\
\sum_{\sigma_1,\sigma_2\in\Sigma_d}\textrm{sgn }\sigma_1\textrm{sgn }\sigma_2a^{\sigma_1(1)\sigma_1(2)\sigma_2(1)\sigma_2(2)}\cdot\cdot\cdot a^{\sigma_1(d-1)\sigma_1(d)\sigma_2(d-1)\sigma_2(d)}.
\end{eqnarray}
\textit{Proof}. By Lemma 5.4.2 we have
\begin{eqnarray}
\textrm{Tr}\Big((-1)^F(\textrm{D}^p\mathcal{A})^\frac{d}{2}\Big)\quad\qquad\qquad\qquad\qquad\qquad\qquad\qquad\qquad\qquad\qquad\nonumber\\
=\textrm{Tr}\Big((-1)^Fa^{i_1j_1i_2j_2}\cdot\cdot\cdot a^{i_{d-1}j_{d-1}i_dj_d}\textrm{D}^p(e_{i_1}^\ast\otimes e_{j_1})\circ\textrm{D}^p(e_{i_2}^\ast\otimes e_{j_2})\Big)\circ\quad\nonumber\\
\cdot\cdot\cdot\circ\textrm{D}^p(e_{i_{d-1}}^\ast\otimes e_{j_{d-1}})\textrm{D}^p(e_{i_d}^\ast\otimes e_{j_d})\qquad\qquad\qquad\qquad\qquad\qquad\qquad\nonumber\\
=\sum_{\substack{i_1,\cdot\cdot\cdot,i_d\\j_1,\cdot\cdot\cdot,j_d}}a^{i_1i_2j_1j_2}\cdot\cdot\cdot a^{i_{d-1}j_{d-1}i_dj_d}\times\qquad\qquad\qquad\qquad\qquad\qquad\qquad\quad\nonumber\\
\textrm{ coefficients of }x^1\cdot\cdot\cdot x^d\textrm{ in }\det(x^ke_{j_k}^\ast\otimes e_{j_k})\quad\qquad\qquad\qquad\qquad\qquad\nonumber\\
=\sum_{\sigma_1,\sigma_2\in\Sigma_d}\textrm{sgn }\sigma_1\textrm{sgn }\sigma_2a^{\sigma_1(1)\sigma_1(2)\sigma_2(1)\sigma_2(2)}\cdot\cdot\cdot a^{\sigma_1(d-1)\sigma_1(d)\sigma_2(d-1)\sigma_2(d)}.\quad
\end{eqnarray}
The proof is complete.\qquad$\square$\\
From now on we shall work with $V=T^\ast M$ to obtain some lemmas concerning the commutation of the covariant derivatives.\\
\textbf{Lemma 5.4.5.} \textit{Suppose} $X_1,\cdot\cdot\cdot,X_m\in C^\infty(U,TM)$, $\mathcal{A}\in C^\infty(U,TM\otimes T^\ast M)$, \textit{then we have the following relation}:
\begin{eqnarray}
\nabla_{X_1}\circ\cdot\cdot\cdot\circ\nabla_{X_m}\circ\textrm{D}^p\mathcal{A}=\textrm{D}^p\mathcal{A}\circ\nabla_{X_1}\circ\cdot\cdot\cdot\circ\nabla_{X_m}\qquad\qquad\qquad\qquad\nonumber\\
\label{eq:com}+\sum_{k=1}^m\sum_{\substack{\sigma\in\Sigma_m\\ \sigma(1)<\cdot\cdot\cdot<\sigma(k)\\ \sigma(k+1)<\cdot\cdot\cdot<\sigma(m)}}\textrm{D}^p\Big(\nabla_{X_{\sigma(1)}}\circ\cdot\cdot\cdot\circ\nabla_{X_{\sigma(k)}}(\mathcal{A})\Big)\circ\nonumber\\
\nabla_{X_{\sigma(k+1)}}\circ\cdot\cdot\cdot\circ\nabla_{X_{\sigma(m)}}.\qquad\qquad\qquad\qquad\qquad\quad
\end{eqnarray}
\textit{Proof}. We will argue by induction. For $m=1$, (\ref{eq:com}) reads
\begin{equation}\label{eq:ini}\nabla_{X_1}\circ\textrm{D}^p\mathcal{A}\omega=\textrm{D}^p\mathcal{A}\circ\nabla_{X_1}\omega+\textrm{D}^p(\nabla_{X_1}\mathcal{A})\omega,\omega\in C^\infty(U,\mathit{\Lambda}^pT^\ast M).\end{equation}
We verify (\ref{eq:ini}) by induction on $p$. For $p=1$, (\ref{eq:ini}) holds by the definition of a covariant derivative. Since the operators $\textrm{D}^p\mathcal{A}$ and $\nabla_{X_1}$ satisfy the axioms of a derivation, one verifies easily that if (\ref{eq:ini}) holds for $\omega_1\in C^\infty(U,\mathit{\Lambda}^p T^\ast M)$ and $\omega_2\in C^\infty(U,\mathit{\Lambda}^q T^\ast M)$, then it also holds for $\omega_1\wedge\omega_2$. This finishes the verification.\\
Suppose $i\in\mathbb{Z}_+$ and the lemma holds for $m\leq i$, by hypothesis we have
\begin{eqnarray}
\nabla_{X_1}\circ\cdot\cdot\cdot\circ\nabla_{X_{i+1}}\qquad\qquad\qquad\qquad\qquad\qquad\qquad\qquad\qquad\qquad\quad\nonumber\\
=\nabla_{X_1}\circ\cdot\cdot\cdot\circ\nabla_{X_i}\textrm{D}^p\mathcal{A}\circ\nabla_{X_{i+1}}+\nabla_{X_1}\circ\cdot\cdot\cdot\circ\nabla_{X_{i}}\circ\textrm{D}^p(\nabla_{X_{i+1}}\mathcal{A})\quad\nonumber\\
=\textrm{D}^p\mathcal{A}\circ\nabla_{X_1}\circ\cdot\cdot\cdot\circ\nabla_{X_{i+1}}+\qquad\qquad\qquad\qquad\qquad\qquad\qquad\qquad\nonumber\\
\bigg(\sum_{k=1}^i\sum_{\substack{\sigma_1\in\Sigma_i\\ \sigma_1(1)<\cdot\cdot\cdot<\sigma_1(k)\\ \sigma_1(k+1)<\cdot\cdot\cdot<\sigma_1(i)}}\textrm{D}^p\Big(\nabla_{X_{\sigma_1(1)}}\circ\cdot\cdot\cdot\circ\nabla_{X_{\sigma_1(k)}}(\mathcal{A})\Big)\quad\qquad\qquad\qquad\nonumber\\
\circ\nabla_{X_{\sigma_1(k+1)}}\circ\cdot\cdot\cdot\circ\nabla_{X_{\sigma_1(i)}}\bigg)\circ\nabla_{X_{i+1}}\qquad\qquad\qquad\qquad\qquad\qquad\qquad\nonumber\\
+\sum_{k=0}^i\sum_{\substack{\sigma_2\in\Sigma_i\\ \sigma_2(1)<\cdot\cdot\cdot<\sigma_2(k)\\ \sigma_2(k+1)<\cdot\cdot\cdot<\sigma_2(i)}}\textrm{D}^p\Big(\nabla_{X_{\sigma_2(1)}}\circ\cdot\cdot\cdot\circ\nabla_{X_{\sigma_2(k)}}(\mathcal{A})\Big)\qquad\qquad\qquad\quad\nonumber\\
\circ\nabla_{X_{\sigma_2(k+1)}}\circ\cdot\cdot\cdot\circ\nabla_{X_{\sigma_2(i)}}\qquad\qquad\qquad\qquad\qquad\qquad\qquad\qquad\qquad\nonumber\\
=\textrm{D}^p\mathcal{A}\circ\nabla_{X_1}\circ\cdot\cdot\cdot\circ\nabla_{X_{i+1}}+\qquad\qquad\qquad\qquad\qquad\qquad\qquad\qquad\nonumber\\
\sum_{k=1}^{i+1}\sum_{\substack{\sigma\in\Sigma_{i+1}\\ \sigma(1)<\cdot\cdot\cdot<\sigma(k)\\ \sigma(k+1)<\cdot\cdot\cdot<\sigma(i+1)}}\textrm{D}^p\Big(\nabla_{X_{\sigma(1)}}\circ\cdot\cdot\cdot\circ\nabla_{X_{\sigma(k)}}(\mathcal{A})\Big)\qquad\qquad\qquad\qquad\nonumber\\
\circ\nabla_{X_{\sigma(k+1)}}\circ\cdot\cdot\cdot\circ\nabla_{X_{\sigma(i+1)}}.\qquad\qquad\qquad\qquad\qquad\qquad\qquad\qquad\qquad
\end{eqnarray}
This proves the lemma.\qquad$\square$\\
The following result is just an anologue of Lemma 5.4.5 under the notation convention which we have explained above.\\
\textbf{Lemma 5.4.6.} \textit{Suppose} $\mathcal{A}\in C^\infty(U,TM\otimes T^\ast M\otimes TM\otimes T^\ast M)$, $X_i,1\leq i\leq m$ \textit{are as above}, \textit{then} (\ref{eq:com}) \textit{holds}.\\
\textit{Proof}. Without loss of generality, we can assume that $\mathcal{A}=\mathcal{B}\otimes\mathcal{C}$ with $\mathcal{B},\mathcal{C}\in C^\infty(U,TM\otimes T^\ast M)$. Then by Lemma 5.4.5, one easily verifies that (\ref{eq:com}) holds for $m=1$. The whole lemma follows by an induction argument on $m$, which is similar with the proof of Lemma 5.4.5, and we shall omit the details.\qquad$\square$\\
Recall that the Riemann curvature tensor $R$ is of type (1,3), therefore it can be identified with a tensor of type (2,2) via the Riemannian metric. Also, it is evident that $R(X,Y)\in C^\infty(U,TM\otimes T^\ast M)$ with $X,Y\in C^\infty(U,TM)$. Thus both $\textrm{D}^pR$ and $\textrm{D}^p\Big(R(X,Y)\Big)$ are well defined.\\
\textbf{Lemma 5.4.7.} \textit{Suppose} $X\in C^\infty(U, TM)$, $X_i,1\leq i\leq m$ \textit{are as above, then}
\begin{eqnarray}\nabla_{X_1}\circ\cdot\cdot\cdot\circ\nabla_{X_m}\circ\nabla_{X}=\nabla_{X}\circ\nabla_{X_1}\circ\cdot\cdot\cdot\circ\nabla_{X_m}+\qquad\nonumber\\
\sum_{j=0}^{m-1}\sum_{\substack{\sigma\in\Sigma_m\\ \sigma(1)<\cdot\cdot\cdot<\sigma(j+1)\\ \sigma(j+2)<\cdot\cdot\cdot<\sigma(m)}}\textrm{D}^p\Big(\nabla_{X_{\sigma(1)}}\circ\cdot\cdot\cdot\circ\nabla_{X_{\sigma(j)}}R(X,X_{\sigma(j+1)})\Big)\nonumber\\
\circ\nabla_{X_{\sigma(j+2)}}\circ\cdot\cdot\cdot\circ\nabla_{X_{\sigma(m)}}+\qquad\qquad\qquad\qquad\qquad\qquad\quad\nonumber\\
\label{eq:cur}\sum_{i=1}^m\nabla_{X_1}\circ\cdot\cdot\cdot\circ\nabla_{X_{i-1}}\circ\nabla_{[X_i,X]}\circ\nabla_{X_{i+1}}\circ\cdot\cdot\cdot\circ\nabla_{X_m},\qquad
\end{eqnarray}
\textit{where} $[\cdot,\cdot]$ \textit{is the Lie bracket}.\\
\textit{Proof}. For $m=1$ and $p=1$, (\ref{eq:cur}) is just the definition of $R$. Then the lemma follows by first applying induction on $p$, then arguing inductively on $m$. Since the proof is similar with that of Lemma 5.4.5, we shall omit the details.\qquad$\square$\\
These commutation lemmas enables us to generalize Lemma 5.4.2 and 5.4.3 to obtain a cancelation lemma which is powerful enough to prove Theorem 5.4.1.\\
\textbf{Lemma 5.4.8. (Patodi's Cancelation Lemma)} \textit{Suppose} $l_1,l_2,l_3,i\in\mathbb{N}$ \textit{satisfying one of the following two conditions}:
\begin{itemize}
\item $l_3>0$, $l_1+2l_2+l_3+2i\leq d$;
\item $l_1+2l_2+l_3+2i<d$.
\end{itemize}
\textit{Let} $\sigma\in\Sigma_{l_1+l_2}$, $j_1,\cdot\cdot\cdot j_{l_3}\in\{1,\cdot\cdot\cdot,d\}$, $\mathcal{A}_{1},\cdot\cdot\cdot,\mathcal{A}_{l_1}\in C^\infty(U,TM\otimes T^\ast M)$, \textit{and} $\mathcal{A}_{l_1+1},\cdot\cdot\cdot,\mathcal{A}_{l_1+l_2}\in C^\infty(U,TM\otimes T^\ast M\otimes TM\otimes T^\ast M)$. \textit{Then}
\begin{equation}\label{eq:main}\sum_{p=0}^d(-1)^p\textrm{tr}_x\bigg(\textrm{D}^p\mathcal{A}_{\sigma(1)}\circ\cdot\cdot\cdot\circ\textrm{D}^p\mathcal{A}_{\sigma(l_1+l_2)}\circ\nabla_{\frac{\partial}{\partial y_{j_1}}}\circ\cdot\cdot\cdot\circ\nabla_{\frac{\partial}{\partial y_{j_{l_3}}}}\Big(u^{i,p}(x,y)\Big)\bigg)(x,x)=0,\end{equation}
\textit{where all the operators act on} $y$.\\
\textit{Proof}. We shall argue by induction on $i$ and $l_3$. For $i=l_3=0$, (\ref{eq:main}) holds by Lemma 5.4.2 and Lemma 5.4.3. Therefore we can suppose the lemma holds for $i=s-1$ for some $s\in\mathbb{Z}_+$, let's consider the case when $i=s$.\\
We first make the following abbreviation:
\begin{equation}\mathcal{A}^p=\textrm{D}^p\mathcal{A}_{\sigma(1)}\circ\cdot\cdot\cdot\circ\textrm{D}^p\mathcal{A}_{\sigma(l_1+l_2)}\circ\nabla_{\frac{\partial}{\partial y_{j_1}}}\circ\cdot\cdot\cdot\circ\nabla_{\frac{\partial}{\partial y_{j_{l_3}}}}.\end{equation}
By (\ref{eq:pre}) and the well-known formula
\begin{equation}\Delta_y=-g^{jk}(y)\nabla_{\frac{\partial}{\partial y^k}}\circ\nabla_{\frac{\partial}{\partial y^j}}+g^{jk}(y)\Gamma_{jk}^\ell(y)\nabla_{\frac{\partial}{\partial y^\ell}}-\textrm{D}^pR_y,\end{equation}
where $R_y$ means the curvature homomorphism acts on $y$, we are able to get an equation which is suitable for induction:
\begin{eqnarray}
\nabla_{r\frac{\partial}{\partial r}}u^{s,p}(x,y)+(s+\frac{r}{4g}\frac{dg}{dr})u^{s,p}(x,y)=\quad\qquad\qquad\qquad\qquad\nonumber\\
\label{eq:ind}\Big(g^{jk}(y)\nabla_{\frac{\partial}{\partial y^k}}\circ\nabla_{\frac{\partial}{\partial y^j}}-g^{jk}(y)\Gamma^\ell_{ij}(y)\nabla_{\frac{\partial}{\partial y^\ell}}+\textrm{D}^pR_y\Big)u^{s-1,p}(x,y).
\end{eqnarray}
We derive from (\ref{eq:pre}) the important fact $\nabla_{r\frac{\partial}{\partial r}}u^{s,p}(x,y)=0$, by first applying $\mathcal{A}^p$ and then taking the trace on each side of (\ref{eq:ind}) we get
\begin{eqnarray}s\sum_{p=0}^d(-1)^p\textrm{tr}_x\bigg(\mathcal{A}^p\Big(u^{s,p}(x,y)\Big)\bigg)(x,x)\qquad\qquad\qquad\qquad\qquad\nonumber\\
=g^{jk}(y)\sum_{p=0}^d(-1)^p\textrm{tr}_x\bigg(\mathcal{A}^p\circ\nabla_{\frac{\partial}{\partial y^k}}\circ\nabla_{\frac{\partial}{\partial y^j}}\Big(u^{s-1,p}(x,y)\Big)\bigg)(x,x)\nonumber\\
-g^{jk}(y)\Gamma_{jk}^\ell(y)\sum_{p=0}^d(-1)^p\textrm{tr}_x\bigg(\mathcal{A}^p\circ\nabla_{\frac{\partial}{\partial y^\ell}}\Big(u^{s-1,p}(x,y)\Big)\bigg)(x,x)\nonumber\\
\label{eq:atr}+\sum_{p=0}^d(-1)^p\textrm{tr}_x\bigg(\mathcal{A}^p\circ\textrm{D}^pR_y\Big(u^{s-1,p}(x,y)\Big)\bigg)(x,x).\qquad\qquad\quad
\end{eqnarray}
Since the right hand side of (\ref{eq:atr}) vanishes by the induction hypothesis, we have
\begin{equation}\label{eq:trace}\sum_{p=0}^d(-1)^p\textrm{tr}_x\bigg(\mathcal{A}^p\Big(\nabla_{r\frac{\partial}{\partial r}}u^{s,p}(x,y)\Big)\bigg)(x,x)=0.\end{equation}
Now suppose $t\in\mathbb{Z}_+$ and the lemma holds for $l_3\leq t-1$. To make Lemma 5.4.7 available, we shall make use of the fact $[\frac{\partial}{\partial y^i},r\frac{\partial}{\partial r}]=\frac{\partial}{\partial y^i}$, which can be verified by direct computation.\\
We finish the proof by first applying $\mathcal{A}^p$ to (\ref{eq:ind}) and then taking the trace. Keeping this goal in mind, we carry out this operation term by term in (\ref{eq:ind}).\\
\textbf{The first term on the left hand side}\\
We derive from (\ref{eq:trace}) the fact that $\nabla_{r\frac{\partial}{\partial r}}\circ\nabla_{\frac{\partial}{\partial y^{j_1}}}\circ\cdot\cdot\cdot\circ\nabla_{\frac{\partial}{\partial y^{j_{l_3}}}}\Big(u^{s,p}(x,y)\Big)=0$, by applying Lemma 5.4.7 and the induction hypothesis to (\ref{eq:trace}), we have
\begin{eqnarray}
\sum_{p=0}^d(-1)^p\textrm{tr}_x\bigg(\mathcal{A}^p\circ\nabla_{r\frac{\partial}{\partial r}}\Big(u^{s,p}(x,y)\Big)\bigg)(x,x)\nonumber\\
=t\sum_{p=0}^d(-1)^p\textrm{tr}_x\bigg(\mathcal{A}^p\Big(u^{s,p}(x,y)\Big)\bigg)(x,x).\qquad
\end{eqnarray}
\textbf{The second term on the left hand side}\\
Direct computation gives
\begin{eqnarray}
\mathcal{A}^p\Big((s+\frac{r}{4g}\frac{dg}{dr})u^{s,p}(x,y)\Big)=\bigg((s+\frac{r}{4g}\frac{dg}{dr})\mathcal{A}^p\Big(u^{s,p}(x,y)\Big)\bigg)+\quad\qquad\nonumber\\
\sum_{k=1}^{l_3}\sum_{\substack{\epsilon\in\Sigma_{l_3}\\ \epsilon(1)<\cdot\cdot\cdot<\epsilon(k)\\ \epsilon(k+1)<\cdot\cdot\cdot<\epsilon(l_3)}}\Big(\nabla_{\frac{\partial}{\partial y^{{j_{\epsilon(1)}}}}}\circ\cdot\cdot\cdot\circ\nabla_{\frac{\partial}{\partial y^{{j_{\epsilon(k)}}}}}(s+\frac{r}{4g}\frac{dg}{dr})\Big)\quad\qquad\qquad\nonumber\\
\textrm{D}^p\mathcal{A}_{\sigma(1)}\circ\cdot\cdot\cdot\circ\textrm{D}^p\mathcal{A}_{\sigma(l_1+l_2)}\circ\nabla_{\frac{\partial}{\partial y^{{j_{\epsilon(k+1)}}}}}\circ\cdot\cdot\cdot\circ\nabla_{\frac{\partial}{\partial y^{{j_{\epsilon(l_3)}}}}}\Big(u^{s,p}(x,y)\Big).
\end{eqnarray}
We then take the trace and use the induction hypothesis to obtain
\begin{eqnarray}\sum_{p=0}^d(-1)^p\textrm{tr}_x\bigg(\mathcal{A}^p\Big((s+\frac{r}{4g}\frac{dg}{dr})u^{s,p}(x,y)\Big)\bigg)(x,x)\nonumber\\
=s\sum_{p=0}^d(-1)^p\textrm{tr}_x\bigg(\mathcal{A}^p\Big(u^{s,p}(x,y)\Big)\bigg)(x,x).\qquad\qquad
\end{eqnarray}
\textbf{Terms on the right hand side}\\
The first two terms on the right hand side can be treated similarly as the second term on the left hand side, i.e., by direct computation and using the induction hypothesis. We easily get
\begin{equation}\sum_{p=0}^d(-1)^p\textrm{tr}_x\bigg(\mathcal{A}^p\circ g^{jk}(y)\nabla_{\frac{\partial}{\partial y^k}}\circ\nabla_{\frac{\partial}{\partial y^j}}\Big(u^{s,p}(x,y)\Big)\bigg)(x,x)=0,\end{equation}
\begin{equation}\sum_{p=0}^d(-1)^p\textrm{tr}_x\bigg(\mathcal{A}^p\circ g^{jk}(y)\Gamma_{jk}^\ell(y)\nabla_{\frac{\partial}{\partial y^\ell}}\Big(u^{s-1,p}(x,y)\Big)\bigg)(x,x)=0.\end{equation}
For the third term, we use Lemma 5.4.6 and the induction hypothesis to get
\begin{equation}\sum_{p=0}^d(-1)^p\textrm{tr}_x\bigg(\mathcal{A}^p\circ\textrm{D}^pR_y\Big(u^{s-1,p}(x,y)\Big)\bigg)(x,x)=0.\end{equation}
Since $t+s>0$, we have $\sum_{p=0}^d(-1)^p\textrm{tr}_x\bigg(\mathcal{A}^p\Big(u^{s,p}(x,y)\Big)\bigg)(x,x)=0$, which finishes the proof.\qquad$\square$\\
\textsc{Proof of Theorem 5.4.1}\\
By Lemma 5.4.8, $\sum_{p=0}^d(-1)^p\textrm{tr}_xu^{i,p}(x,x)=0$ for $i<\frac{d}{2}$, we then use (\ref{eq:ind}) to compute $\sum_{p=0}^d(-1)^p\textrm{tr}_xu^{\frac{d}{2},p}(x,x)$. Apply Lemma 5.4.8 twice to get
\begin{eqnarray}
\sum_{p=0}^d\textrm{tr}_xu^{\frac{d}{2},p}(x,x)=\frac{1}{n}\sum_{p=0}^d\textrm{tr}_x\bigg(\textrm{D}^pR_y\Big(u^{\frac{d}{2}-1,p}(x,y)\Big)\bigg)(x,x)\qquad\qquad\qquad\nonumber\\
=\frac{1}{n(n-1)}g^{jk}(x)\sum_{p=0}^d(-1)^p\textrm{tr}_x\bigg(\textrm{D}^pR_y\circ\nabla_{\frac{\partial}{\partial y^k}}\circ\nabla_{\frac{\partial}{\partial y^j}}\Big(u^{\frac{d}{2}-2,p}(x,y)\Big)\bigg)(x,x)\nonumber\\
=\frac{1}{n(n-1)}g^{jk}(x)\Gamma_{jk}^\ell(x)\sum_{p=0}^d(-1)^p\textrm{tr}_x\Big(\nabla_{\frac{\partial}{\partial y^\ell}}u^{\frac{d}{2}-2,p}(x,y)\Big)(x,x)\qquad\qquad\quad\nonumber\\
=\frac{1}{n(n-1)}\sum_{p=0}^d(-1)^p\textrm{tr}_x\bigg(\textrm{D}^pR_y\circ\textrm{D}^pR_y\Big(u^{\frac{d}{2}-2,p}(x,y)\Big)\bigg)(x,x).\qquad\qquad\quad
\end{eqnarray}
By Lemma 5.4.8, the first two terms on the right hand side again vanish, therefore we get
\begin{eqnarray}
\sum_{p=0}^d(-1)^p\textrm{tr}_xu^{\frac{d}{2},p}(x,x)\qquad\qquad\qquad\qquad\qquad\qquad\qquad\qquad\nonumber\\
=\frac{1}{n(n-1)}\sum_{p=0}^d(-1)^p\textrm{tr}_x\bigg(\textrm{D}^pR_y\circ\textrm{D}^pR_y\Big(u^{\frac{d}{2}-2,p}(x,y)\Big)\bigg)(x,x).
\end{eqnarray}
Proceeding like this we finally get
\begin{equation}\sum_{p=0}^d(-1)^p\textrm{tr}_xu^{\frac{d}{2},p}(x,x)=\frac{1}{n!}\sum_{p=0}^d(-1)^p\bigg((\textrm{D}^pR_y)^{\frac{d}{2}}\Big(u^{0,p}(x,y)\Big)\bigg)(x,x).\end{equation}
Combine Lemma 5.4.4 and the first Bianchi indentity, the proof is complete.\qquad$\square$\\
\textsc{Remark}\\
In \cite{ms}, McKean and Singer gave their formula for the index of the Dirac operators (see also \cite{bgv}):
\begin{equation}\label{eq:msfd}\textrm{ind}_A(\mathsf{D})=\textrm{Tr}\Big((-1)^Fe^{-t\mathsf{D}^2}\Big),\end{equation}
where $\textrm{ind}_A$ denotes the analytic index, and $\mathsf{D}$ denotes a Dirac operator.\\
Substitute $\mathsf{D}$ by $d+\delta$ in (\ref{eq:msfd}), the Gauss-Bonnet-Chern theorem then follows as a corollary of the local index theorem.\\
Actually, based on (\ref{eq:msfd}), the McKean-Singer conjecture can be raised for all Dirac operators, and a proof of the generalized McKean-Singer conjecture yields the local Atiyah-Singer index theorem. We refer the reader to \cite{bgv} and \cite{eg} for more details.

\begin{flushleft}
Department of Mathematics, NanJing Normal University, 210046, China\\
\textit{Email address}: \textsf{yinlee1004@msn.com}
\end{flushleft}

\end{document}